\documentclass{article}
\title{Shock Formation of 3D Euler-Poisson System for Electron Fluid with Steady Ion Background}
\date{}
\author{Yiya Qiu\footnote{School of Mathematical Sciences, University of Science and Technology of China, Hefei, Anhui, 230026, PR China, Email:qq171579@mail.ustc.edu.cn}, Lifeng Zhao\footnote{School of Mathematical Sciences, University of Science and Technology of China, Hefei, Anhui, 230026, PR China, Email:zhaolf@ustc.edu.cn}}

\makeatletter
\@addtoreset{equation}{section}
\makeatother

\usepackage{xcolor}

\usepackage{amsmath}
\usepackage{amssymb}
\usepackage{amsthm}
\usepackage{mathtools}
\usepackage[a4paper,left=3cm,right=3cm,top=3cm,bottom=3cm]{geometry}
\usepackage{subeqnarray}
\usepackage{cases}
\usepackage{setspace}
\usepackage{bbm} 
\usepackage{enumerate}

\newtheorem{thm}{Theorem}[section]
\newtheorem{lem}[thm]{Lemma}
\newtheorem{cor}[thm]{Corollary}
\newtheorem{prop}[thm]{Proposition}

\newtheorem{rek}[thm]{Remark}

\newcommand\al{\alpha} 
\newcommand\be{\beta}  
\newcommand\ga{\gamma}  
\newcommand\ve{\varepsilon}  
\newcommand\la{\lambda}    
\newcommand\de{\delta}   
\newcommand\ka{\kappa} 

\newcommand\ze{\zeta}

\newcommand\si{\sigma}

\def\XX{{\scriptstyle \mathcal{X} }}

\newcommand*{\supp}{\ensuremath{\mathrm{supp\,}}}

\renewcommand*{\tilde}{\widetilde}

\renewcommand*{\bar}{\overline}

\newcommand{\xcal}{{\text{x}}}

\newcommand{\tcal}{{\text{t}}}
\newcommand{\ck}{\check}

\newcommand\les{\lesssim} 
\newcommand\f{\frac}  
\newcommand\p{\partial}  

\newcommand\td{\tilde} 
\newcommand\wt{\widehat} 

\allowdisplaybreaks
\begin{document}
\maketitle
\begin{abstract}
We prove shock formation for 3D Euler-Poisson system for electron fluid in plasma. The shock solution we construct is of large initial data, also compactly supported during the lifespan. In addition, the blowup time and location can be computed explicitly.
\end{abstract}

\section{Introduction}

We consider the 3D repulsive Euler-Poisson system for electron fluid in plasma
\begin{equation}\label{EP}
\begin{cases}
n_em_e(\p_\tcal u+u\cdot\nabla_\xcal u)+\nabla_\xcal p_e=-n_e e\nabla_\xcal\phi,\\
\p_\tcal n_e+\nabla_\xcal\cdot(n_e u)=0,\\
\Delta_\xcal\phi=4\pi e(n_+-n_e).
\end{cases}
\end{equation}
Here, we denote charge of electron by e, mass of electron by $m_e$,  electron density by $n_e:\mathbb{R}^3\times\mathbb{R}\to\mathbb{R}$,  velocity by $u:\mathbb{R}^3\times\mathbb{R}\to\mathbb{R}^3$ and  pressure by $p_e=p(n_e)=\f{1}{\ga}n_e^\ga$ with $\ga>1$. The electric potential $\phi:\mathbb{R}^3\times\mathbb{R}\to\mathbb{R}$ depends on both electrons fluid $n_e$ and ion background $n_+$. In addition, it satisfies $\phi(\xcal)\to0$ as $|\xcal|\to \infty.$

The Euler-Poisson system arises in various physical backgrounds involving compressible fluids interacting with a self-consistent potential. In stellar dynamics, Euler-Poisson system is used to describe the evolution of self-gravitational gaseous stars in which the potential results from the attractive gravitational interaction. In plasma physics, Euler-Poisson system is a simplified version of the so called  ``two-fluid'' model which describes plasma dynamics for ion and electron fluids. In such case, the potential comes from repulsive Coulomb interaction. Because the ratio of the electron mass and the ion mass is very small, the heavier ions are treated as motionless in some cases and the density of ions can be regarded as a time independent function. 

There have been a lot of works devoted to the study of Euler-Poisson system.
For two dimensional Euler-Poisson system \eqref{EP}, the global existence of smooth irrotational flows with small initial amplitude was proved in \cite{Ionescu2013}, \cite{Jang2012}, \cite{Jang2014} and \cite{Li2014}. In three dimensional case, Guo constructed the global irrotational flows with small velocity for the electron fluid in \cite{Guo1998}.
There are also many results for other types of Euler-Poisson system in plasma. In \cite{Guo2011}, Guo and Pausader established the global existence and uniqueness for ion fluid with Boltzman statistics, but whether the shock waves can develop in 2D remains open. What is more, for full Euler-Poisson system concerning plasma namely two-fluid model, the global existence for small solutions was established  by Guo, Ionescu and Pausader in \cite{Guo2016} and \cite{Guo2014} via the method of space-time resonance.

The situation is different for attractive Euler-Poisson system in gaseous stars. In the star-formation model, the range of index $\ga$ in the pressure law $p(\ga)=\rho^\ga$ is essential. An important question is the dynamical stability of Lane-Emden stars solutions. In \cite{Rein2003}, Rein showed the existence of global weak solutions by variational methods and obtained stability of Lane-Emden stars for $\ga>\f43$. On the contrary, unstability was proved by Jang for $\ga=\f65$  in \cite{Jang2008}, and for $\f65<\ga<\f43$ in \cite{Jang2014a}. Global expanding stars are also well studied in recent years. Hadzi\'c and Jang discovered in \cite{Hadzic2017} that the expanding solutions are asymptotically stable when $\ga=\f43$. In addition, Liu \cite{Liu2019} studied the viscous effect in expanding configuration. The gravitational collapse is another fascinating issue. The key is to find the self-similar solutions, and they have been constructed in \cite{Deng2003}, \cite{Goldreich1980} and \cite{Makino1992}. Besides, in \cite{Guo2021}, gravitational collapse was shown by Guo, Hadzi\'c and Jang based on pressureless model. Furthermore, a very recent work \cite{Guo2021a} showed the existence of smooth
radially symmetric self-similar solutions to the gravitational Euler-Poisson system when $1<\ga<\f43$.

However, the blowup mechnism for Euler-Poisson system needs further investigation. In fact, there has been a lot of remarkable results for compressible Euler in recent years. For 3D Euler system, Sideris argued by contradiction and proved in \cite{Sideris1985} that $C^1$ regular solutions of Euler have a finite lifespan. The first proof of shock formation was given by Chirstodoulou for relativistic fluid in \cite{Chris2007} and non-relativistic fluid by Christodoulou and Miao in \cite{Christodoulou2014},  both of them concerns irrotational flows. Luk and Speck then proved for shock formation for 2D fluid flows with vorticity in \cite{Luk2018}. Buckmaster, Shkoller and Vical constructed the stable asymptotically self-similar type blowup solutions to 2D and 3D compressible Euler in \cite{Buckmaster2020}, \cite{Buckmaster2019}. They also analyzed the shock formation in \cite{Buckmaster2020a} for the 3d non-isentropic Euler system, in which sounds waves interact with entropy waves to produce vorticity. 
Unstable shock solutions are constructed by Buckmaster and Iyer for 2D Euler system in \cite{Buckmaster2019a} and 
 Buckmaster,  Drivas,  Shkoller and Vicol
 proved that a discontinuous shock instantaneously
develops after the pre-shock in \cite{Buckmaster2021}.
  In addition, Merle, Raphael, Rodnianski and Szeftel constructed smooth spherical blowup solutions for compressible Euler and Navier-Stokes system by analyzing phase portrait in \cite{Merle2019} and \cite{Merle2019a}. For Euler-Poisson system, the study of shock formation is of great interests. On one hand, it has been shown in \cite{Guo1998}, \cite{Guo2016a}, \cite{Ionescu2013}, \cite{Li2014} that smooth solutions with small amplitude to 
Euler–Poisson system for electrons persist forever with no shock formation. On the other hand, Makino and Perthame proved in \cite{Makino1990} that smooth solutions with gravitational force blow up in finite time if the initial data are spherically symmetric and have a compact support. In addition, Guo and Tahvildar-Zadeh proved in \cite{Guo1999} that shock waves do develop for ``large'' perturbations of the constant state equilibrium of $n_e\equiv0$, $u\equiv0$. Moreover, in \cite{Wang2014}, Wang showed that for a larger class of initial data, with no matter repulsive forces or attractive, the finite time blowup take place. The assumption on the compact support was removed by Li and Wang in \cite{Li2018}.

In this article we study the shock formation for Euler-Poisson system (\ref{EP}) for large initial data. Our main results can be stated roughly as:

\begin{thm}
	The system (\ref{EP}) admits solution
	 $(u,n_e)$ developing shock in finite time. More precisely, the minimum negative slopes of $u_1$ and $n_e$ go to $-\infty$ in finite time, while themselves remain bounded.
\end{thm}
The idea is mainly inspired by the pioneering work \cite{Buckmaster2019} for 3D compressible Euler system. After Riemann transformation, we get the system (\ref{Riewzu}) about $(w,z,u_\nu)$ see Section 2, which is regarded as a perturbation of 3D Burgers equation
\begin{equation}\label{3db}
\p_t u+u\cdot\nabla u=0.
\end{equation}
It is widely acknowledged that (\ref{3db}) admits a family of self-similar solutions, which are the prototypes of various shock formation problems and universally appear in the different circumstances of fluid dynamics, see \cite{Eggers2009}.
Here we still expect that the singular solutions we construct to (\ref{Riewzu}) is of self-similar type, and behave just like that of (\ref{3db}), namely $\bar W$ defined in (\ref{barW}). In order to verify the assertion rigorously we need to perform our argument in self-similar variables $(y,s)$ and inverstigate the difference $\td W=W-\bar W$, where $W$ is the singular solution we construct in self-similar variable $(y,s)$.

However, due to the presence of the electric potential, (\ref{EP}) does not enjoy finite speed of propagation as compressible Euler system in \cite{Buckmaster2019}. Indeed, the absence of finite speed of propagation is usually caused by nonlocal term in various fluid dynamics, such as Biot-Savart law in incompressible flow, or nonlocal drift. As a result, when estimating nonlinear terms or external forces one has to take the nonlocal effect into consideration. In \cite{Yang2020}, Yang gave a good example analyzing external forces with nonlocal terms.

Fortunately, when it comes to Euler-Poisson system, the situation changes, because the integrand in nonlocal term is just the electron density $n_e$. On one hand, the transport of $n_e$ relies heavily on the velocity, which stays bounded in $L^\infty$(see Lemma 4.7) during the whole lifespan. That means if $n_{e,0}$ is supported in compact set $\XX_0$ at initial time, then for any $\tcal>-\f{2}{1+\al}$, $n_e(\xcal,\tcal)$ is supported in compact set $\psi(\XX_0)$, where $\psi$ is Lagrangian trajectory generated by velocity field. On the other hand, the velocity field $u$ is also compactly supported for all time, because nontrivial velocity cannot be attached to vacuum. As a result, we can perform estimates just in the compact set $\mathcal X(s)$ in self-similar variables.

In our construction, due to the various symmetries of the system, some unstable modes need to be considered. They are parameters $\ka$, $\tau$, $\xi_1$, $\xi_2$, $\xi_3$, correspongding to Galilean transform, time translation and spatial translation. In order to give precise description of blowup solutions, we impose some constrains on $W(0,s)$, $\nabla W(0,s)$ and $\nabla^2 W(0,s)$. In fact, these parameters are required to satisfy proper ODEs and we will prove these ODEs are solvable, see Proposition 4.7. But we should point out that  the authors in \cite{Buckmaster2019}, except the above five parameters, introduced another five parameters $n_2$, $n_3$, $\phi_{22}$, $\phi_{23}$ and $\phi_{33}$ concerning the rotation symmetry and second fundamental form of wave front, forcing the amount $\ck\nabla W(0,s)= \ck\nabla^2W(0,s)\equiv0$. In this paper, in order to avoid the over-repeating parts of \cite{Buckmaster2019} and focus on the differences, we impose the the ansatz  that $u$ and $n_{e,0}$ are even about $\xcal_2$ and $\xcal_3$ in (\ref{even}) and $\ck\nabla^2(u_{1,0}+\f{n_{e,0}^\al}{\al})|_{\tcal=-\f{2\ve}{1+\al}}=0$ in (\ref{ckna}) so that these five modes are excluded. As a consequence, there are more unstable modes than that of \cite{Buckmaster2019} for our blowup solutions. Of course, if we involve these five modes in our argument, similar conclusions can be obtained systematically as in \cite{Buckmaster2019}.

As the end of the introduction, we claim that the similar shock formation result for 2D Euler-Poisson system can also be obtained by following the work \cite{Buckmaster2020}.

The remainder of the paper is organized as follows: In Section 2, we perform a series of transforms from sound speed form (\ref{usigma}) to self-similar form (\ref{WZU}), then we apply higher order derivatives and get (\ref{Devo}). Also, similar process is taken to the perturbation $\td W=W-\bar W$. In section 3, we make some assumptions on the initial data and give the main theorem, both of which are stated in self-similar variables. In section 4, we make bootstrap assumptions for no more than 4 order derivatives, and give the $\dot H^k$ estimate in Proposition 4.2 for $k\geq 18$. Also, we analyze the modulation variables. In section 5, we estimate the bounds of transport terms and external forces under the bootstrap assumptions. In section 6, we close the bootstrap assumptions made in Proposition 4.1. In section 7, we prove the energy estimate Proposition 4.2. In section 8, we state the main theorem in physical variables.

\textbf{Notations.} 

For a vector $A=(a_1,a_2,a_3)\in\mathbb{R}^3$, we always use $\ck a$ to denote anyone of $a_2$ and $a_3$, as well as $|\ck a|=(a_2^2+a_3^2)^\f12$, and for muti-index $\ga=(\ga_1,\ck \ga) $ is similar. We let the subscript $i,j,k\in\{1,2,3\}$, and Greek subscript $\mu, \nu\in\{2,3\}$.

The inequality $A\les B$ means $A\leq C B$, where $C$ is a constant depending on $\al$ and $\ka_0$, but independent of $M$, $\ve$, $l$ and $L$. Concerning the various constants we suppose $0<\ve\ll l\leq\f{1}{10}<1\ll M\ll L=\ve^{-\f{1}{10}}$, and for any $\de_1,$ $\de_2>0$, there exists $\de_3>0$, such that $\ve^{\de_1}M^{\de_2}<\ve^{\de_3}$, where $\de_3<\de_1.$

\section{The Reformulation of the Main Problem}

We assume that $n_{e,0}=n_e(\xcal,-\f{2}{1+\al}\ve)$ and $u_0=u(\xcal,-\f{2}{1+\al}\ve)$ are supported in the set 

\begin{equation}\label{sptx}
	\XX=\{|\xcal_1|\leq\ve^\f12, |\ck\xcal|\leq \ve^{\f16} \}
\end{equation}
and $n_{e,0}>0$ and assume that 
\begin{subequations}\label{even}
\begin{align}
u_0(\xcal_1,\xcal_2,\xcal_3)=u_0(\xcal_1,-\xcal_2,\xcal_3)=u_0(\xcal_1,\xcal_2,-\xcal_3)\\
n_{e,0}(\xcal_1,\xcal_2,\xcal_3)=n_{e,0}(\xcal_1,-\xcal_2,\xcal_3)=n_{e,0}(\xcal_1,\xcal_2,-\xcal_3)
\end{align}
\end{subequations}

Also, for background $n_+$ we further assume that  $n_+$ is supported in 
\begin{equation}\label{sptx+}
	\XX_+=\{|\xcal|\leq1\}
\end{equation}
and lies in Sobolev space $W^{4,\infty}$ with 
\begin{equation}\label{n+}
	\|n_+\|_{W^{4,\infty}}\leq1.
\end{equation}
Besides, the total charges satisfy neutrality  $$\int_{\XX_+}\big(n_+-n_e\big)=0.$$

The vorticity $\omega_e=\nabla\times u$, so the specific vorticity $\zeta_e=\f{\omega_e}{n_e}$ satisfies 
\begin{equation}\label{zeta}
	\p_t\zeta_e+(u_e\cdot \nabla_\xcal )\zeta_e-(\zeta_e\cdot\nabla_\xcal)u_e=0.
\end{equation}

If we let sound speed $\sigma_e=\f{1}{\al}n_e^\al$, where $\al=\f{\ga-1}{2},$ $m_e=e=1$, then we have the system about $u$ and $\sigma_e$
\begin{equation}\label{usigma}
	\begin{cases}
		\p_\tcal u+(u\cdot\nabla_\xcal)u+\al\si_e\nabla_\xcal\si_e=-\nabla_\xcal\phi,\\
		\p_\tcal\si_e+(u\cdot\nabla_\xcal)\si_e+\al\si_e\nabla_\xcal\cdot u=0.
	\end{cases}
\end{equation}

\subsection{A Series of Transforms}
In order to investigate the structure of equation in depth, various  transforms performed on the origin system (\ref{usigma}) are needed. Next we will see the blow up solution we construct are close to the stable self-similar solution of 3D Burgers, while the blowup location is related to the amplitude of velocity field. To this end, we need three steps to obtain the final form of the system we study. The key point to figure out the different forms of electric potential produced by charges.

\textbf{Step 1. Translation} 

At first, we define $$t=\f{1+\al}{2}\tcal$$ and assume the blowup take place at the space-time point $(\xcal,t)=(\xi(T_*), T_*)$, where $\xi(t)$ is a function depending only on $t$.  In fact, $\xi(t)$ is used to track the location of the most steepened point of the gradient, which we assume initially satisfies $\xi(-\ve)=0$. Let $$x=\xcal-\xi(t),\ \ \td{u}(x,t)=u(\xcal,t),\ \ \td{\si}(x,t)=\si_e(\xcal,t),$$ $$\td{n}_e(x,t)=n_e(\xcal,t),\ \td{n}_+(x,t)=n_+(\xcal,t),\ \ \td\zeta(x,t)=\zeta_e(\xcal,t).$$
Then the potential function $\phi$ is translation invariant, that is
\begin{align*}
\phi(\xcal)=\int\f{1}{|\bar{\xcal}|}(n_+-n_e)(\xcal-\bar{\xcal})d\bar{\xcal}=
\int\f{1}{|\bar{x}+\xi|}(n_+-n_e)(x+\xi-\bar{\xcal})d\bar{x}=\td{\phi}(x).
\end{align*}
Therefore, the system (\ref{usigma}) is reformulated as
\begin{equation}\label{xiusgm}
\begin{cases}
\f{1+\al}{2}\p_t\td{u}+((\td{u}-\dot\xi)\cdot\nabla)\td{u}+\al\td{\si}\nabla\td{\si}=-\nabla\td{\phi},\\
\f{1+\al}{2}\p_t\td{\si}+((\td{u}-\dot\xi)\cdot\nabla)\td{\si}+\al\td{\si}\nabla\cdot\td{u}=0.
\end{cases}
\end{equation}

\textbf{Step 2. Riemann Transform}

We define the constants
$$\be_1=\f{1}{\al+1},\ \ \be_2=\f{1-\al}{1+\al},\ \ \be_3=\f{\al}{1+\al},$$  and introduce the Riemann variables
$$w=\td{u}_1+\td{\si},\ \  z=\td{u}_1-\td{\si},\ \ \mu=2,3.$$ 
Then the system (\ref{xiusgm}) is rewritten as
\begin{subequations}\label{Riewzu}
\begin{align}
\p_t w+(w+\be_2z-2\be_1\dot\xi_1)\p_1w+2\be_1(\td{u}_\nu-\dot\xi_\nu)\p_\nu w=-2\be_3\td{\si}\p_\nu \td{u}_\nu-2\be_1\p_1\td{\phi},\\
\p_t z+(\be_2w+z-2\be_1\dot\xi_1)\p_1z+2\be_1(\td{u}_\nu-\dot\xi_\nu)\p_\nu z=2\be_3\td{\si}\p_\nu \td{u}_\nu-2\be_1\p_1\td{\phi},\\
\p_t\td{u}_\mu+(\be_1w+\be_1z-2\be_1\dot\xi_1)\p_1\td{u}_\mu+2\be_1(\td{u}_\nu-\dot\xi_\nu)=2\be_3\td{\si}\p_\mu\td{\si}-2\be_1\p_\mu\td{\phi}.
\end{align}
\end{subequations}

\textbf{Step 3. Self-Similar Transform}

Self-similar transform is widely used in the study of blowup phenomenon in PDE theory. By zooming in the area near the blowup point corresponding to blowup speed, we can give the precise description for the profile, asymptotic behavior, stability properties.

In this problem, we make the ansatz that the shock accumulation of 3D  Euler-Poisson system appears like 3D Burgers equation, so the blowup speed inherits from the shock formation of Burgers equation.

We define the self-similar variables
\begin{equation*}
s=-\log(\tau(t)-t),\ \  y_1=e^{\f{3}{2}s}x_1,\ \  y_\nu=e^{\f{s}{2}}x_\nu,\ \ \be_\tau=\f{1}{1-\dot\tau}
\end{equation*}
and set $$w(x,t)=e^{-\f{s}{2}}W(y,s)+\ka(t),\ \ z(x,t)=Z(y,s),\ \ \td{u}_i(x,t)=U_i(y,s),$$ as well as  $$\td{\si}(x,t)=S(y,s),\ \ \td n_e(x,t)=R_e(y,s),\ \ \td n_+(x,t)=R_+(y,s),\ \ \td\zeta(x,t)=\Omega(y,s),$$ where we assume initially $$\tau(-\ve)=0,\ \  \ka(-\ve)=\ka_0>1.$$ 
Meanwhile the electric potential $\Phi$ in the self-similar variables $(y,s)$ is
\begin{equation}\label{poten}
\Phi(y,s)=\int\f{(R_+-R_e)(y-z)e^{-\f{5}{2}s}}{(e^{-3s}|z_1|^2+e^{-s}|\check{z}|^2)^{\f{1}{2}}}dz.
\end{equation} 
Now by (\ref{sptx+}) and (\ref{n+}), we have
\begin{equation}\label{sptX+}
\supp R_+\subset \mathcal X_+=\{|y_1|\leq e^{\f32s},\ |\ck y|\leq e^{\f s2} \}
\end{equation}
\begin{equation}\label{pR+}
\|\p^\ga R_+\|_{L^\infty}\leq e^{-\f{3\ga_1}{2}s-\f{|\ck\ga|}{2}s},\ \ \ |\ga|\leq4
\end{equation}
Therefore, we obtain from (\ref{Riewzu}) that
\begin{subequations}\label{WZU}
\begin{align}
(\p_s-\f{1}{2})W+(g_W+\f{3}{2}y_1)\p_{y_1}W+(h^\mu+\f{1}{2}y_\mu)\p_{y_\mu}W=&F_W-\be_\tau e^{-\f{s}{2}}\dot\ka,\label{evoW}   \\
\p_sZ+(g_Z+\f{3}{2}y_1)\p_{y_1}Z+(h^\mu+\f{1}{2}y_\mu)\p_{y_\mu}Z=&F_Z,\\
\p_s U_i+(g_U+\f{3}{2}y_1)\p_{y_1}U_i+(h^\mu+\f{1}{2}y_\mu)\p_{y_\mu}U_i=&F_{U_i}, \label{DevoU}\\
\p_s S+(g_U+\f{3}{2}y_1)\p_{y_1}S+(h^\mu+\f{1}{2}y_\mu)\p_\mu S=&F_S,\label{DevoS}
\end{align}
\end{subequations}
where the transport velocity fields are
\begin{subequations}\label{tspt}
\begin{align}
g_W=\be_\tau W+G_W, G_W=\be_\tau e^{\f{s}{2}}(\ka+\be_2Z-2\be_1\dot\xi_1),\label{GW}\\
g_Z=\be_2\be_\tau W+G_Z, G_Z=\be_\tau e^{\f{s}{2}}(\be_2\ka+Z-2\be_1\dot\xi_1),\\
g_U=\be_1\be_\tau W+G_U, G_U=\be_\tau e^{\f{s}{2}}(\be_1\ka+\be_1Z-2\be_1\dot\xi_1),\\
h^\mu=\be_\tau e^{-\f{s}{2}}  (2\be_1U_\mu-2\be_1\dot\xi_\mu)\label{h}
\end{align}
\end{subequations}
and the external forces are
\begin{subequations}
\begin{align}
F_W=-2\be_3\be_\tau S(\p_\nu U_\nu)-2\be_1\be_\tau e^s \p_1\Phi,\label{FW1}\\
F_Z=2\be_3\be_\tau e^{-\f{s}{2}}S(\p_\nu U_\nu)-2\be_1\be\tau e^{\f{s}{2}}\p_1\Phi,\\
F_{U_\nu}=-2\be_3\be_\tau e^{-\f{s}{2}}S\p_\nu S-2\be_1\be_\tau e^{-\f{s}{2}}\p_\nu\Phi,\label{FUnu1}\\
F_{U_1}=-2\be_3\be_\tau e^{\f{s}{2}}S\p_1 S-2\be_1\be_\tau e^{\f{s}{2}}\p_1\Phi,\\
F_S=-2\be_\tau\be_3S(e^{\f{s}{2}}\p_{y_1}U_1+e^{-\f{s}{2}}\p_{y_\nu}U_{\nu}).
\end{align}
\end{subequations}

\subsection{The Evolution of Higher Order Derivatives}
In the study of quasilinear system, it is necessary to consider the evolution of higher order derivatives.
 Apply $\p^\ga$ to the equations of $(W,Z,U_i)$, with $|\ga|\geq1$, and $\ga=(\ga_1,\ck\ga)$, we get
\begin{subequations}\label{Devo}
\begin{align}
\left( \p_s+\f{3\ga_1+\ga_2+\ga_3-1}{2}+\be_\tau(1+\ga_1\mathbbm{1}_{\ga_1\geq2})\p_1W \right) \p^\ga W+(\mathcal{V}_W\cdot\nabla)\p^\ga W=F^{(\ga)}_W\label{Devo1},\\
\left( \p_s+\f{3\ga_1+\ga_2+\ga_3}{2}+\be_2\be_\tau\ga_1\p_1W  \right) \p^\ga Z+(\mathcal{V}_Z\cdot\nabla)\p^\ga Z=F^{(\ga)}_Z,\label{Devo2}\\
\left( \p_s+\f{3\ga_1+\ga_2+\ga_3}{2}+\be_1\be_\tau\ga_1\p_1W  \right) \p^\ga U_\nu+(\mathcal{V}_Z\cdot\nabla)\p^\ga U_\nu=F^{(\ga)}_{U_\nu},\label{Devo3}
\end{align}
\end{subequations}where 
\begin{subequations}
\begin{align}
 \mathcal{V}_W=(g_W+\f{3}{2}y_1, h^2+\f{1}{2}y_2, h^3+\f{1}{2}y_3),\\
 \mathcal{V}_Z=(g_Z+\f{3}{2}y_1, h^2+\f{1}{2}y_2, h^3+\f{1}{2}y_3),\\
 \mathcal{V}_U=(g_U+\f{3}{2}y_1, h^2+\f{1}{2}y_2, h^3+\f{1}{2}y_3).
\end{align}
\end{subequations}
And the external forces are given by 
\begin{subequations}
\begin{align}
F^{(\ga)}_W=&\p^\ga F_W-\sum_{0\leq\be<\ga}\binom\ga\be\left(\p^{\ga-\be}G_W\p_1\p^\be W+\p^{\ga-\be}h^\mu\p_\mu\p^\be W \right)\notag \\
&-\mathbbm{1}_{|\ga|\geq2}\be_\tau\sum_{\substack{|\be|=|\ga|-1\\ \ga_1=\be_1}}\p^{\ga-\be}W\p^\be \p_1W-\mathbbm{1}_{|\ga|\geq3}\be_\tau\sum_{\be\leq\ga-2}\binom\ga\be \p^{\ga-\be}W\p^\be \p_1W, \label{FW2}\\
F^{(\ga)}_Z=&\p^\ga F_Z-\sum_{0\leq\be<\ga}\binom\ga\be\left(\p^{\ga-\be}G_Z\p_1\p^\be Z+\p^{\ga-\be}h^\mu\p_\mu\p^\be Z \right)\notag\\
&-\be_2\be_\tau\sum_{\substack{|\be|=|\ga|-1\\ \ga_1=\be_1}}\p^{\ga-\be}W\p^\be \p_1Z-\mathbbm{1}_{|\ga|\geq2}\be_2\be_\tau\sum_{\be\leq\ga-2}\binom\ga\be \p^{\ga-\be}W\p^\be \p_1Z, \label{FZ2}\\
F^{(\ga)}_{U_i}=&\p^\ga F_{U_i}-\sum_{0\leq\be<\ga}\binom\ga\be\left(\p^{\ga-\be}G_{U}\p_1\p^\be {U_i}+\p^{\ga-\be}h^\mu\p_\mu\p^\be {U_i} \right)\notag\\
&-\be_1\be_\tau\sum_{\substack{|\be|=|\ga|-1\\ \ga_1=\be_1}}\p^{\ga-\be}W\p^\be \p_1U_i-\mathbbm{1}_{|\ga|\geq2}\be_1\be_\tau\sum_{\be\leq\ga-2}\binom\ga\be \p^{\ga-\be}W\p^\be \p_1U_i.
\end{align}
\end{subequations}

\subsection{The Perturbation around the Burgers Profile}
In order to describe the asymptotic profile for $W(s,y)$, a solution of \eqref{evoW}, we rely heavily on the the stable self-simlar solutions to 3D Burgers equation
\begin{equation}\label{3dbur}
-\f12\bar{W}+(\f32y_1+\bar{W})\p_1\bar{W}+\f12y_\mu\p_{\mu}\bar{W}=0.
\end{equation}
An example of such solutions is
\begin{equation}\label{barW}
\bar{W}(y)=(1+|\ck{y}|^2)^{\f12}W^*((1+|\ck{y}^2|)^{-\f32}y_1),
\end{equation}
where $W^*$ solves 1D Burgers self-similar equation
$$-\f12 W^*+(\f32y+W^*)\p_yW^*=0.$$
Moreover, Buckmaster, Shkoller and Vicol in \cite{Buckmaster2019}  gave the asymptotic  estimates of $\bar{W}$ by means of the weighed function $\eta(y)=(1+|y_1|^2+|\ck y|^6)$
\begin{equation}\label{pbarW}
\|\eta^{-\f16}\bar{W}\|_{L^\infty}\leq1,\ 
\|\eta^{\f13}\p_1\bar{W}\|_{L^\infty}\leq 1,\ 
\|\ck\nabla\bar{W}\|_{L^\infty}\leq 1,\ 
\|\eta^{\f13}\p_1\nabla\bar{W}\|_{L^\infty}\leq1,\ 
\|\eta^{\f16}\ck\nabla^2\bar{W}\|_{L^\infty}\leq1.
\end{equation}

In fact, There are more solutions to (\ref{3dbur}) than just (\ref{barW}). In \cite{Chris2007}, Christodoulou defined the non-degeneracy condition to classify somehow genuine 3d shock formation in the quasilinear system. In \cite{Buckmaster2019}, the authors introduced the equivalent notion called ``genericity'', by checking whether the Hessian matrix $\p_1\nabla^2\bar{W}$ is positive definite. The family of stable generic solutions of (\ref{3dbur}) can be indexed by a 3-tensor $\mathcal{A}_\al$, where $\al$ is a multi-index with $|\al|=3$. In particular, we have
\begin{prop}\emph{(\cite{Buckmaster2019})}
Let $\mathcal{A}$ be symmetric a 3-tensor such that $\mathcal{A}_{1jk}=\mathcal{M}_{jk}$ with $\mathcal{M}$ is a positive definite symmetric matrix. Then there exists a $C^\infty$ solution $\bar{W}_\mathcal{A}$ to 
$$-\f12\bar{W}_\mathcal{A}+\left(\f32y_1+\bar{W}_\mathcal{A} \right)+\f{\ck y}{2}\cdot\ck\nabla\bar{W}_\mathcal{A}=0,$$
which satisfies
\begin{itemize}
\item
$\bar{W}_\mathcal{A}(0)=0, \p_1\bar{W}_\mathcal{A}(0)=-1, \ck\nabla\bar{W}_\mathcal{A}(0)=0$,
\item 
$\p^\ga \bar{W}_\mathcal{A}(0)$ for $|\al|$ even,
\item
$\p^\ga \bar{W}_\mathcal{A}(0)=\mathcal{A}_\al$ for $|\al|=3.$
\end{itemize}
\end{prop}

In view of the above proposition,  (\ref{barW}) is indeed a generic solution because 
\begin{equation*}
\nabla^2\p_1\bar{W}=\textbf{diag}\{6,2,2\}>0.
\end{equation*}

Now we ignore the differences between each $\bar{W}_\mathcal{A}$ and
compute the evolution of perturbation. Let $\td{W}=W-\bar{W}$, it follows from (\ref{evoW}) and (\ref{3dbur}) that
\begin{align*}
\p_s\td{W}+(\be_\tau\p_1\bar{W}-\f12)\td{W}+(\mathcal{V}_W\cdot\nabla)\td{W}=&F_W-e^{-\f s2}\be_\tau \dot{\ka}+((\be_\tau-1)\bar{W}-G_W)\p_1\bar{W}-h^\mu\p_\mu\bar{W}
=:\td{F}_W.
\end{align*}
Furthermore, the higher derivatives satisfies
\begin{align}\label{DtdW}
\left(\p_s+\f{3\ga_1+\ga_2+\ga_3-1}{2}+\be_\tau(\p_1\bar{W}+\ga_	1\p_1W ) \right)\p^\ga\td{W}+(\mathcal{V}_W\cdot\nabla)\p^\ga\td{W}=\td{F}^{(\ga)}_W,
\end{align}
where $\td{F}^{(\ga)}_W$ is 
\begin{align}
\td{F}^{(\ga)}_W=&\p^\ga \td{F}_W-\sum_{0\leq\be<\ga}\binom\ga\be\left(\p^{\ga-\be}G_W\p_1\p^\be\td{W}+\p^{\ga-\be}h^\mu\p_\mu\p^\be\td{W}+\be_\tau\p^{\ga-\be}(\p_1\bar{W})\p^\be\td{W} \right)\notag\\
&-\be\tau\sum_{\substack{1\leq|\be|\leq|\ga|-2\\ \be\leq\ga}}\binom\ga\be\p^{\ga-\be}W\p_1\p^\be\td{W}-\be\tau\sum_{\substack{|\be|=|\ga|-1\\ \be_1=\ga_1}}\binom\ga\be\p^{\ga-\be}W\p_1\p^\be\td{W}.
\label{FW3}
\end{align}

\section{Main Theorem}
\subsection{Assumptions on the Initial Data }
For the modulation variables, we have assumed in physical variables
\begin{equation}\label{inimod}
\ka|_{t=-\ve}=\ka_0>0,\ \ \tau|_{t=-\ve}=0,\ \ \xi|_{t=-\ve}=0
\end{equation}in Section 2.
Next we mainly state the assumptions on the initial data in self-similar variables. Note that (\ref{sptx}) implies that the initial data for $(W,Z,U_{\nu})$ and $R, S$ are supported in the set 
\begin{equation} \label{inispt}
\mathcal{X}_0=\left\{|y_1|\leq \ve^{-1}, |\ck y|\leq \ve^{-\f13}  \right\}.
\end{equation}
Then in order to match $W$ and $\bar{W}$ at $y=0$, we assume
\begin{equation}
W(0,-\log \ve)=0,\ \ \p_1W(0,-\log\ve)=-1,\ \ \ck\nabla W(0,-\log\ve)=0,\ \ \nabla^2W(0,-\log\ve)=0.
\end{equation}
We define $\td{W}=W-\bar{W}$. Then we get
\begin{equation}
\td{W}(y, -\log\ve)=W(y,-\log\ve)-\bar{W}(y).
\end{equation}

Since $\bar{W}$ does not decay at infinity and $W(y,-\log\ve)$ is supported compactly, $\td W=W-\bar W$ can only be expected to stay close to zero in a smaller region $\{|y|\leq L\ll \ve^{-\f13}\}$.  Therefore,  the assumptions on the initial data and the bootstrap assumptions for $\td{W}$ in the following should be restricted in the region $\{|y|\leq L\}$. Furthermore, if $y$ is very close to the origin, $|y|\leq l$, then it is convenient to use Taylor expansion to give the estimates. Therefore, the assumptions on $\td{W}$ are made in three different pieces:  
\begin{equation*}
	|y|\leq l, \ \ l<|y|\leq L, \ \ |y|>L,
\end{equation*}
where \begin{equation}\label{lL}
	l=M^{-\f{1}{200}},\ \ L=\ve^{-\f{1}{10}}.
\end{equation}

We assume the initial data of $W,Z,U_\nu$, specific vorticity $\Omega$, and their $\dot H^k$ norm as following 
\begin{enumerate}[i).]
	\item For $|y|\leq L$, we assume that 
	\begin{subequations}\label{initdW}
		\begin{align}
			\eta^{-\f16}(y)|\td{W}(y,-\log\ve)|&\leq\ve^\f14,\\
			\eta^{\f13}(y)|\p_1\td{W}(y,-\log\ve)|&\leq\ve^\f14,\\
			|\ck\nabla\td{W}(y,-\log\ve)|&\leq\ve^\f14.
	\end{align}\end{subequations}
\item For $|y|\leq l\ll1$, we assume
\begin{subequations}
	\begin{align}
		|\p^\ga\td{W}(y,-\log\ve)|\leq M\ve^{\f{1}{18}}, &\  |y|\leq l,\ \ |\ga|=4\label{initdW4},\\
		|\p^\ga\td{W}(0,-\log\ve)|\leq M\ve^{\f{1}{18}}, &\  |\ga|=3\label{initdW3}.
	\end{align}
\end{subequations}
\item Because we have (\ref{initdW}), (\ref{pbarW}) and $W=\bar{W}+\td{W}$ when $|y|\leq L$,
so we just make assumptions for $\p^\ga W$ with $|\ga|\leq1$ in the region $\{y\geq L\}\cap\mathcal{X}_0$.
\begin{subequations}
	\begin{align}
		\eta^{-\f16}(y)|W(y,-\log\ve)|\leq&1+\ve^{\f{1}{10}},\label{iniW}\\
		\eta^{\f13}(y)|\p_1W(y,-\log\ve)|\leq&1+\ve^{\f{1}{10}},\label{iniW1}\\
		|\ck\nabla W(y,-\log\ve)|\leq&1+\ve^{\f{1}{10}}\label{iniWck}.
	\end{align}
\end{subequations}
\item For the second derivatives of $W$, for all $y\in\mathcal{X}_0$, we have
\begin{subequations}
	\begin{align}
		\eta^{\f13}(y)|\p^\ga W(y,-\log\ve)|\leq1,&\ \ \ga_1\geq1,\ |\ga|=2,\label{iniW2}\\
		\eta^{\f16}(y)|\ck\nabla^2W(y,-\log\ve)|\leq1\label{iniW2ck}&.
	\end{align}
\end{subequations}

\item
For $Z$, we assume 
\begin{equation}\label{iniZ}
|\p^\ga Z(y,-\log\ve)|\leq
\begin{cases}
\ve^{\f32}, &\ga_1\geq1\ |\ga|=1,2\\
\ve,  &\ga_1=0\ |\ck\ga|=0,1,2
\end{cases}
\end{equation}

\item
For $U_\nu$, we assume
\begin{equation}\label{iniUnu}
|\p^\ga U_\nu(y,-\log\ve)|\leq
\begin{cases}
\ve^{\f32}, &\ga_1=1\ |\ck\ga|=0\\
\ve,  &\ga_1=0\ |\ck\ga|=0,1,2
\end{cases}
\end{equation}

\item
For the initial specific vorticity, we assume that
\begin{equation}\label{iniOme}
\|\Omega_0\|\leq1.
\end{equation}

\item
At last, we assume the homogenous Sobolev norm for initial data of $W,\ Z,\ U_\nu$ and $R_+$ satisfy
\begin{equation}\label{iniHk}
\ve\|W(\cdot,-\log\ve)\|^2_{\dot{H}^k}+\|Z(\cdot,-\log\ve)\|^2_{\dot{H}^k}+\|U_\nu(\cdot,-\log\ve)\|^2_{\dot{H}^k}\leq\ve
\end{equation}
and 
\begin{equation}\label{R+Hk}
\|R_+\|^2_{\dot H^k}  \leq \ve
\end{equation}
for all $k\geq 18$.
\end{enumerate}	

\subsection{The Statement of Main Theorem}
\begin{thm}
\emph{(Shock formation in self-similar variables)}
Let $\nu=2,3$, $\ga>1$, $\al=\f{\ga-1}{2},$ $\ka_0=\ka_0(\al)>1.$ We set the initial time is $s=-\log\ve$ and consider the system \emph{(\ref{WZU}a)-(\ref{WZU}c)} of $(W,Z,U_\nu)$. Assume the initial data $(W,Z,U_\nu)|_{s=-\log\ve}$ satisfy the assumptions \emph{(\ref{inispt})}-\emph{(\ref{iniHk})}. Then there exist $M(\al,\ka_0)\gg1$ and $0<\ve(\al,\ka_0,M)\ll1$ such that the system \emph{(\ref{WZU}a)-(\ref{WZU}c)} admits a unique solution $(W,Z,U_\nu)\in C([-\log\ve,\infty);H^k)$ defined in $\mathcal{X}(s)$ \emph{(\ref{sptX})}, with $k\geq18$.  Furthermore, it holds that
\begin{equation}\label{wzuhk}
\|W(\cdot,s)\|^2_{\dot{H}^k}+e^s\|Z(\cdot,s)\|^2_{\dot{H}^k}+\|U_\nu(\cdot,s)\|^2_{\dot{H}^k}\leq \la^{-k}e^{-s-\log\ve}+(1-e^{-s-\log\ve})M^{4k}.
\end{equation}

Besides, $W$ behave similarly to the stable self-similar profile $\bar{W}$ of 3D Burgers. More precisely, if $\td{W}=W-\bar{W}$, then
$\td{W}$ satisfies \emph{(\ref{ptdWW})-(\ref{ptdWW3})}, as well as $\p^\ga\td{W}(0,s)=0$ for $|\ga|\geq2$. In addition, $W$ satisfies \emph{(\ref{pW}), (\ref{pW3})}. In fact, the limit function $\bar{W}_\mathcal{A}(y)=\lim_{s\to\infty}W(y,s)$ are uniquely determined by $\p^\al W(0,s)$ for $|\al|=3$.

The size of $Z$ and $U_\nu$ stay in $O(\ve)$ for $s\geq-\log\ve$. In particular, they satisfy \emph{(\ref{pZ}), (\ref{pUnu}), (\ref{pZ3}) and  (\ref{pUnu3})}.

The sound speed $S(y,s)$ satisfies
$$\|S(\cdot,s)-\f{\ka_0}{2}\|_{L^\infty}\leq \ve^{\f18}.$$

\end{thm}

\begin{rek}
Compared with \emph{\cite{Buckmaster2019}}, not only $(W,Z,U_\nu)|_{-\log\ve}$ but also $R|_{s=-\log\ve}$ and $S|_{s=-\log\ve}=\f{1}{\al}R^\al|_{s=-\log\ve}$ are required to be compactly supported. Because only if electron density is compactly supported, can system still keep the finite speed of propagation, then the electric field effect caused by charges in compact set is under control in this framework.
\end{rek}
\begin{rek}
Theorem 3.1 holds for all $\ga>1$, as long as taking appropriate $\ka_0$ in \emph{(\ref{ka0})} and $\ve$ sufficiently small. Unlike the cases where the range of $\ga$ plays an essential role in Euler-Poisson system, Theorem 3.1 suggests that the dynamics of shock formation dominates the pressure and repulsive force due to the electric potential, as long as the gradient of initial data is steep enough.
\end{rek}

\begin{rek}
In fact the evolution of specific vorticity $\zeta$ has similar form to that in Euler system, since $\nabla\times\nabla\phi=0$. Combined with the similar bounds in the bootstrap assumption, there holds the same conclusion obtained in \emph{\cite{Buckmaster2019}}. That is there exists $C>0$, such that $\f1C|\Omega(y_0)|^2<\Omega(\Psi^{y_0}_U(s),s)<C|\Omega(y_0)|^2$ for all $s\geq-\log\ve$, where $\Psi^{y_0}(s)$ is the Lagrangian trajectory generated by velocity $\mathcal V_U.$
\end{rek}
\proof[Proof of Theorem 3.1]
The condition of support set  (\ref{sptX}) will be proved in Lemma 4.6. (\ref{wzuhk}) will be proved in Proposition 4.4. The bounds of $\p^\ga \td{W}$, $\p^\ga W$, $\p^\ga Z$ and $\p^\ga U_\nu$ will be proved in Section 7.

The proof of the convergence $\lim_{s\to\infty} W(y,s)=\bar{W}_\mathcal{A}$ is very similar to that in \cite{Buckmaster2019}. The only difference here is due to the electric potential. However, it plays an ignorable role in external forces.

\qed

\section{Bootstrap Assumptions}
\subsection{Assumptions on Derivatives}

In this section, we make bootstrap bounds on derivatives of $W$,$\td W$, $Z$ and $U_\nu$, as well as the support. Note that the bootstrap constants we choose in the following differ from order to order, from variable to variable, even vary by several orders of magnitude. There is a profound connection between these constants and the structure of equation. The appropriate choice of these constants is essential and provides technical convenient in the proof.

We should point out that the following three subsections and the external force estimates as their corollary in the next section are parallel, so the cross-reference is safety as long as the bootstrap constant can be narrowed in the argument.

\begin{prop}\emph{(Main bootstrap assumption)}
We have the following bounds for $\p^\ga W$, $\p^\ga \td{W}$, $\p^\ga Z$ and $\p^\ga U_\nu$:
\begin{enumerate}[i).]
\item \textbf{$W$ bootstrap.} 

First, we postulate the following derivative estimates of $W$:
\begin{equation}\label{pW}
|\p^\ga W|\leq 
\begin{cases}
2\eta^{\f16}(y), & |\ga|=0\\
2\eta^{-\f13}(y),  & \ga_1=1, |\check{\ga}|=0\\
2, &\ga_1=0, |\check{\ga}|=1\\
2M^{\f{1+|\ck\ga|}{3}}\eta^{-\f13}(y), & \ga_1\geq1, |\ga|=2\\
2M\eta^{-\f16}(y), & \ga_1=0, |\check{\ga}|=2
\end{cases}
\end{equation}
Next, we assume that for $|y|\leq L$, the following bounds hold
\begin{subequations}\label{ptdWW}
	\begin{align}
	|\td{W}(y,s)|\leq \ve^{\f{1}{20}}\eta^{\f16}(y),\ \ \ |y|\leq L\label{ptdW}\\
	|\p_1\td{W}(y,s)|\leq \ve^{\f{1}{20}}\eta^{-\f13}(y),\ \ \ |y|\leq L \label{ptdW1}\\
	|\ck{\nabla}\td{W}(y,s)|\leq \ve^{\f{1}{20}}, \ \ \ |y|\leq L\label{ptdWck}.
	\end{align}
\end{subequations}

Furthermore, for $|y|\leq l$ we assume that
\begin{subequations}\label{ptdWW3}
	\begin{align}
	|\p^\ga \td{W}|\leq &\f23M\ve^{\f{1}{18}}l^{4-|\ga|},                      \ &|\ga|=3 \label{ptdW3}\\
	|\p^\ga \td{W}|\leq& \f23M\ve^{\f{1}{18}}, &|\ga|=4 \label{ptdW4}
	\end{align}
\end{subequations}
while for $y=0$, we have 
\begin{equation}\label{ptdWW03}
|\p^\ga \td{W}(0,s)|\leq \ve^{\f{1}{10}},\ \ \ \ \ \ \  |\ga|=3
\end{equation}

\item \textbf{$Z$ bootstrap.} 

We postulate the following derivative estimates of $Z$:
\begin{equation} \label{pZ}
|\p^\ga Z|\leq M
\begin{cases}
e^{-\f32 s}, & \ga_1\geq1\ \ |\ga|=1,2\\
\ve^\f12 , & |\ga|=0\\
\ve^\f12e^{-\f12s} &\ga_1=0,\ |\ck\ga|=1\\
e^{-s} &\ga_1=0,\ |\ck\ga|=2
\end{cases}
\end{equation}

\item \textbf{$U$ bootstrap.} 

We postulate the following derivative estimates of $U_\nu$:
\begin{equation}\label{pUnu}
|\p^\ga U_\nu|\leq M
\begin{cases}
e^{-\f32s}, &\ga_1=1,\ |\check{\ga}|=0,\\
\ve^{\f12},&|\ga|=0\\
\ve^{\f{1}{2}}e^{-\f{1}{2}s}, &\ga_1=0,\ |\ck\ga|=1\\
e^{-s}, &\ga_1=0,\ |\ck\ga|=2
\end{cases}
\end{equation}
\end{enumerate}
\end{prop}

In fact, from the second inequality of \emph{(\ref{pW})} and \emph{(\ref{ptdW1})}, we can get 
\begin{equation}\label{pWlbd}
|\p_1W|\leq 1+\ve^{\f{1}{20}}.
\end{equation}
Indeed, when $|y|\geq L$, $2\eta^{-\f13}\ll1$. When $l\leq |y|\leq L$, by (\ref{pbarW}) and (\ref{ptdW1}), $|\p_1W|\leq (1+\ve^{\f{1}{20}})\eta^{-\f13}\leq1+\ve^{\f{1}{20}}$.

The $L^\infty$ type bounds are not enough to close the bootstrap argument because the higher order derivatives are present in nonlinear terms. So we need to perform $\dot H^k$ estimates and use Sobolev interpolation to compensate the loss of  derivates.

\begin{prop}\label{hk}
For integer $k\geq18$, we have
\begin{subequations}\label{WZUhk}
\begin{align}
\|Z\|^2_{\dot{H}^k}+\|U_\nu\|^2_{\dot{H}^k}\leq 2\la^{-k}e^{-s}+e^{-s}(1-e^{-s}\ve^{-1})M^{4k},\\
\|W\|^2_{\dot{H}^k}\leq 2\la^{-k}e^{-s}\ve^{-1}+(1-e^{-s}\ve^{-1})M^{4k},
\end{align}
\end{subequations}
\end{prop}
\proof We will prove Proposition 4.2 in Section 7.

\qed

Proposition \ref{hk} yields the higher order derivatives estimates as follows:
\begin{cor}There hold that
\begin{equation} \label{pW3}
|\p^\ga W(y,s)|\les
\begin{cases}
e^{\f{2s}{2k-7}}\eta^{-\f13}(y), &\ga_1\neq0, |\ga|=3,4\\
e^{\f{s}{2k-7}}\eta^{-\f16}(y), &\ga_1=0, |\ga|=3,4.
\end{cases}
\end{equation}

\begin{equation} \label{pZ3}
|\p^\ga Z(y,s)|\les
\begin{cases}
e^{-(\f32-\f{3}{2k-7})s}, &\ga_1\geq1, |\ga|=3\\
e^{-(1-\f{|\ga|-1}{2k-7})s}, &|\ga|=3,4,5,
\end{cases}
\end{equation}
\begin{equation} \label{pUnu3}
|\p^\ga U_\nu(y,s)|\les
\begin{cases}
e^{-(\f32-\f{2|\ga|-1}{2k-5})s},   &\ga_1\geq1, |\ga|=2,3\\
e^{-(1-\f{|\ga|-1}{2k-7})s}, & |\ga|=3,4,5,
\end{cases}
\end{equation}

\end{cor}

\proof 
The proof is just repeatition of that in \cite{Buckmaster2019} and we omit the details. 

\qed

Since $U_1=\f12(e^{-\f s2}W+\ka+Z)$, $S=\f12(e^{-\f s2}W+\ka-Z)$ and $S=(\al R_e)^\f{1}{\al},$
the estimates of $U_1$, $S$ and $R_e$ can be seen as a corollary of Proposition 4.1.
\begin{cor}
For $y\in\mathcal{X}(s)$ we have
\begin{equation} \label{pU1S}
|\p^\ga U_1|+|\p^\ga S|\les
\begin{cases}
1, & |\ga|=0\\
e^{-\f s2}\eta^{-\f13}(y), &  \ga=(1,0,0)\\
e^{-\f s2}, & \ga_1=0, |\ck\ga|=1\\
Me^{-\f s2}\eta^{-\f13}(y), & \ga_1\geq1, |\ga|=2\\
Me^{-\f s2}\eta^{-\f16}(y), &\ga_1=0, |\ck\ga|=2\\
e^{(-\f12+\f{3}{2k-7})s}\eta^{-\f13}(y), &\ga_1\neq0, |\ga|=3,4\\
e^{(-\f12+\f{3}{2k-7})s}\eta^{-\f16}(y), &\ga_1=0, |\ck\ga|=3,4
\end{cases}
\end{equation}

\begin{equation} \label{pR}
|\p^\ga R_e|\les
\begin{cases}
1, & |\ga|=0\\
M^\f1\al e^{-\f s2}\eta^{-\f13}(y), & \ga_1\geq1, |\ga|=1,2\\
 e^{-\f s2}, & \ga_1=0, |\ck\ga|=1\\
M^\f1\al e^{-\f s2}\eta^{-\f16}(y), &\ga_1=0, |\ck\ga|=2\\
M^\f1\al e^{(-\f12+\f{3}{2k-7})s}\eta^{-\f13}(y), &\ga_1\neq0, |\ga|=3,4\\
M^\f1\al e^{(-\f12+\f{3}{2k-7})s}\eta^{-\f16}(y), &\ga_1=0, |\ga|=3,4
\end{cases}
\end{equation}
while for $|y|\leq l$ and $|\ga|=4$ we have 
$$|\p^\ga U_1|+|\p^\ga S|\les e^{-\f s2}.$$
\end{cor}
\proof 
Note that 
$$|\p^\ga U_1|+|\p^\ga S|\les e^{-\f s2}|\p^\ga W+\ka\mathbbm{1}_{|\ga|=0}|+|\p^\ga Z|$$
and $$|\p^\ga R_e|\les |S^{\f1\al-1}||\p^\ga S|\les M^{\f1\al-1}|\p^\ga S|.$$
Then (\ref{pR}) and (\ref{pU1S}) follow from Proposition 4.1.

\qed

\subsection{The Evolution of Support Set}
We assume that $(W,Z,U_\nu,R_e)$ have compact support
\begin{equation}\label{sptX}
\mathcal{X}(s)=\left\{|y_1|\leq2\ve^\f12e^{\f32s}, |\ck y|\leq2\ve^\f16e^{\f s2}   \right\},\ \ \ s\geq -\log\ve.
\end{equation}
Then if $y\in\mathcal{X}(s)$, the weighed function $\eta(y)$ satisfies
\begin{equation}\label{eta}
\eta^{\f13}(y)\leq4\ve^\f13e^s.
\end{equation}
Since the evolution of $n_e$ does not involve nonlocal term, $n_e$ is expected to be compactly supported for all time once it is assumed at the initial time. 

\begin{lem}For any $s\geq-\log\ve$, 
$S$ is supported in the set $\mathcal{X}(s)$, as long as $n_{e,0}$ is supported in $\XX$.
\end{lem}
\proof
From $x=\xcal-\xi(t)$, $\td\si(x,t)=\si_e(\xcal,\tcal)=\f{1}{\al}n^\al_e(\xcal,\tcal)$, we know that if $n_{e,0}$ is supported in $\XX$, then
$\td\si(x_0,-\ve)$ is supported in $\XX=\{|x_1|\leq\ve^\f12,|\ck x|\leq \ve^\f16\}.$
Suppose $\varphi(x_0,t)$ is the Lagrangian trajectories generated by velocity field $\f{2}{1+\al}(\td u-\dot{\xi})$ and $\varphi(x_0,-\ve)=x_0\in\XX$, from (\ref{xiusgm}) $$\f{1+\al}{2}\p_t \td\si+((\td u-\dot\xi)\cdot\nabla)\td\si=-\al\td\si\nabla\cdot\td u,$$
we get that
$$\td\si(\varphi (x_0,t),t)=\td\si(x_0, 0)\exp\left(-\f{2\al}{1+\al}\int_{-\ve}^{t}(\nabla\cdot \td u)\circ\varphi(x_0,\de)d\de\right),$$
which implies that if $\td\si(\varphi(x_0,t),t)\neq0$, then $\td\si(x_0, 0)\neq0.$ Namely, for any $x_*$ with $\td \si(x_*,t)\neq0$, there exists $x_0\in\XX$, such that $x_*=\varphi(x_0,t).$
But (\ref{modu}) and the bootstrap bounds of $U_1$ and $U_\nu$ imply that $$|x_*-x_0|=|\varphi(x_0,t)-\varphi (x_0,0)|<\|\td u-\dot{\xi}\|_{L^\infty}T_*<\ve^\f12,$$ which means $x_0\notin\{|x|\leq1\}= \XX$. Thus we reach a contradiction.

\qed

The definitions of $\mathcal X_+$ in (\ref{sptX+}) and $\mathcal X(s)$ in (\ref{sptX}) imply that $\mathcal X(s)\subset\mathcal X_+$ for any $s\geq -\log\ve$.

\subsection{Analysis of Modulation Variables}

\qed

In this section we will impose the following constraints on derivatives of $W$ at $y=0$ by choosing appropriate $\ka$, $\tau$ and $\xi$. In order to do this, we first make some preparations. Plugging
\begin{equation}\label{conW}
W(0,s)=0,\ \p_1W(0,s)=-1,\ \check{\nabla}W(0,s)=0,\ \nabla^2W(0,s)=0
\end{equation}
into 
 (\ref{Devo1}) and (\ref{FW2}), then evaluating at $y=0$ we get
\begin{equation}\label{dotka}
-G^0_W=F^0_W-\be_\tau e^{-\f s2}\dot{\ka},\ \ \ \ \ \ |\ga|=0
\end{equation}and
\begin{equation}\label{FG}
F^{(\ga),0}_W=\p^\ga F^0_W+\p^\ga G^0_W,\ \ \ \ \ \ |\ga|=1,2.
\end{equation}

By the evolution of $\p_1 W$ at $y=0$, we have 
$$-(1-\be_\tau)=\p_1F^0_W+\p_1 G^0_W,$$ which implies that
\begin{equation}\label{dottau}
\dot{\tau}=\f{1}{\be_\tau}(\p_1F^0_W+\p_1G^0_W).
\end{equation}

Next, by the evolution of $\p_1\nabla W$ at $y=0$, it holds
$$F^{(2,0,0)}_W=\p_{11}F_W-\p_{11}G_W\p_1W-\p_{11}h^\mu\p_\mu W-2\p_1G_W\p_{11}W-2\p_1h^\mu\p_{\mu}W.$$
For $i=1,2,3$, we have
\begin{equation*}
G^0_W\p_{1i1}W^0+h^{\mu,0}\p_{1i\mu}W^0=\p_{1i}F^0_W+\p_{1i}G^0_W.
\end{equation*}

Due to (\ref{ptdWW03}), we have $\mathcal{H}^0=\p_1\nabla^2W^0=\p_1\nabla^2\bar{W}^0+\p_1\nabla^2\td{W}^0=\textbf{diag}\{6,2,2\}+\ve^{\f{1}{10}}$,
so $\mathcal{H}^0$ is invertible and
$$(\mathcal{H}^0)^{-1}\leq 1.$$
Thus,
\begin{align}
|G^0_W|+|h^{\mu,0}|\leq& (\mathcal{H}^0)^{-1}(|\p_1\nabla G^0_W|+|\p_1\nabla F^0_W|)\notag\\
\leq& |\p_1\nabla G^0_W|+|\p_1\nabla F^0_W|.\label{G0h0}
\end{align}
 
\begin{lem}For $s\geq-\log\ve$, there holds
\begin{equation}\label{Gh0}
|G^0_W(s)|+|h^{\mu,0}(s)|\leq Me^{-s}.
\end{equation}
\end{lem}
\proof
By (\ref{G0h0}), it suffices to estimate $\p_1\nabla G^0_W$ and $\p_1\nabla F^0_W$. Indeed, (\ref{pG}), (\ref{pZ}) and (\ref{pFWZ}) imply 
$$|\p_1\nabla G^0_W|\leq e^\f s2|\p_1\nabla Z|\leq Me^{-s}$$ and 
$$|\p_1\nabla F^0_W|\leq e^{-s}\eta^{-\f16+\f{2|\ga|-1}{3(2k-5)}}<e^{-s}.$$

\qed

This lemma plays a vital important role in our argument. On one hand, (\ref{Gh0}) holds under the constrains (\ref{conW}), which are guaranteed by the choosing  modulation variables appropriately. On the other hand, (\ref{Gh0}) implies the transport effect in the evolution of $\p^\ga W$ is ignorable, compared with the linear damping term, as we will see in the following sections.

Next we will perform the bootstrap argument to the modulation variables, by using bootstrap bounds assumed in Section 4.1 and their direct corollary: bounds of external forces, whose proof will be postponed in the next section.

\begin{prop}Under the bootstrap assumptions, the modulation variables satisfy the estimates
\begin{equation}\label{modu}
|\dot{\tau}|\leq 2Me^{-s},\ \ \ |\dot\ka|\leq M,\ \ \ |\dot\xi|\leq M^\f12,\ \ \ |T_*|\leq 3M\ve^2.
\end{equation}
Moreover, it holds that
\begin{equation}\label{modu2}
\f12\ka_0\leq|\ka|\leq 2\ka_0\leq M,\ \ \ |\xi|\leq M \ve
\end{equation}
and
\begin{equation}\label{betau}
|1-\be_\tau|=\f{|\dot{\tau}|}{1-\tau}\leq2Me^{-s}\leq2M\ve.
\end{equation}
\end{prop}
\proof

We estimates $\dot\xi$ first. Due to (\ref{GW}) and (\ref{h}), we have
$$
2\be_1\dot\xi_1=-\f{1}{\be_\tau} G^0_We^{-\f s2}+\ka+\be_2 Z^0
$$
and 
$$2\be_1\dot{\xi}_\nu=-\f{e^{\f s2}}{\be_\tau}h^{\mu,0}+2\be_1U_\mu .$$
Hence by (\ref{Gh0}), (\ref{pZ}) and (\ref{pUnu}), it follows that 
 $$|\dot\xi_1|+|\dot\xi_\mu|\leq Me^{-\f {3s}{2}}+\ka_0+M\ve^\f12+Me^{-\f s2}\leq 2\ka_0\leq \f12M^\f12$$ by taking $M$ large enough, which improves the bound in (\ref{modu}).

Next we turn to the estimate of $\dot\ka$. By (\ref{pU1S}), (\ref{pR}) and (\ref{pUnu}), we have
\begin{align*}
|F^0_W|\leq &2\be_3\be_\tau |S^0(\p_\nu U^0_\nu)|+2\be_1\be_\tau e^s|(\p_1\Phi)^0|\\
\leq& M\ve^\f12e^{-\f s2}+e^{-\f s2}\\
\leq &2e^{-\f s2},
\end{align*}
which, together with (\ref{dotka}) and (\ref{Gh0}), yields
$$|\dot\ka|\leq \be_\tau e^{\f s2}(|F^0_W|+|G^0_W|)\leq (1+M\ve)(Me^{-\f s2}+2)\leq \f12M,$$ improving the bound in (\ref{modu}).

Now we are left for the estimate of $\dot\tau$.  (\ref{pZ}), (\ref{dottau}), (\ref{pG}) and (\ref{pFWZ}) imply that 
$$|\dot\tau|\leq |\p_1F^0_W|+|\p_1G^0_W|\les e^{-s}+Me^{-s}\leq  \f32Me^{-s}.$$
Thus the bound in (\ref{modu}) is improved.

Now it remains to estimate $T_*$. From the relation $-\log(\tau(t)-t)=s$ and $\tau(T_*)=T_*$, we have
$$\int_{-\ve}^{T_*}(1-\dot\tau(t))dt=\ve,$$ which implies that 
\begin{equation*}
(T_*+\ve)\leq\ve+\|\dot\tau\|_{L^\infty}(T_*+\ve).
\end{equation*}
Therefore, we get from the estimate of $\dot{\tau}$ in \eqref{modu} that
\begin{equation*}
T_*\leq 2M\ve T_*+2M\ve^2.
\end{equation*}
Thus $T_*\leq 2M\ve^2$ as long as $\ve$ is taken small enough, which improves the bound in (\ref{modu}).

Finally, we are left with  (\ref{modu2}) and (\ref{betau}). Combining initial data (\ref{inimod}), integration $\dot{\ka}$ and $\dot{\xi}$ from $-\ve$ to $T_*$ in (\ref{modu}) gives (\ref{modu2}). And (\ref{betau}) follows from (\ref{modu2}) directly.

\qed

\section{Estimates for Transport Fields and External Forces}
In this section, we compute the bounds of transport fields and external forces. 
\subsection{Transport Term}
\begin{lem}
For $s\geq-\log\ve$ and $y\in\mathcal X(s)$, we have 

\begin{equation}\label{pG}
|\p^\ga G_W|+|\p^\ga G_Z|+|\p^\ga G_U|\leq e^{\f s2}|\p^\ga Z|,\ \ \ |\ga|\geq1
\end{equation}
and 
\begin{equation}\label{ph}
|\p^\ga h^\mu|\leq e^{-\f s2}|\p^\ga U_\mu|, \ \ \ |\ga|\geq1.
\end{equation}
Meanwhile, when $|\ga|=0$, a sharper bound for $G_W$ can be obtaind for any $y\in\mathcal X(s)$
\begin{equation}\label{G}
|G_W|\les M( e^{-\f s2}+|y_1|e^{-s}+\ve^\f12|\ck y|).
\end{equation}
\end{lem}
\proof

(\ref{pG}) and (\ref{ph}) follow from (\ref{tspt}). It remains to prove \eqref{G}. By mean value theorem and (\ref{pG}), 
\begin{align*}
|G_W|\les & |G_W^0|+|y_1|e^{\f s2}\|\p_1 Z\|_{L^\infty}+|\ck y|e^{\f s2}\|\ck \nabla Z\|_{L^\infty}\\
\les&  M(e^{-s}+|y_1|e^{-s}+\ve^{\f12}|\ck y|),
\end{align*}
where in the last inequality we have used (\ref{Gh0}) and (\ref{pZ}).
\qed

\subsection{External Forces with Electric Potential}
Now we estimate the derivatives of electric potential. In the proof the compact support of $R_e$ and $R_+$ is of vital importance, because the Poisson kernel in 3D is not integrable in whole space. 
\begin{lem}
For any $y\in \mathcal X(s)$, and $|\ga|\leq 4$, we have
\begin{equation} \label{pPhi1}
|\p^\ga\p_1\Phi(y)|\les e^{-\f32s}\left(\|\p^\ga R_e\|_{L^\infty}+\|\p^\ga R_+\|_{L^\infty}\right)
\end{equation}and 
\begin{equation}\label{pPhick}
|\p^\ga\p_\nu\Phi(y)|\les e^{-\f s2}(\|\p^\ga R_e\|_{L^\infty}+\|\p^\ga R_+\|_{L^\infty}).
\end{equation}
\end{lem}
\proof
By (\ref{poten}), we have
$$\p_1\Phi=\int_{\mathbb{R}^3}\f{-e^{-\f {5s}{2}}e^{-3s}z_1(R_+-R_e)(y-z)}{(e^{-3s}|z_1|^2+e^{-s}|\ck z|^2)^{\f32}}dz$$ and 
$$\p_\nu\Phi=\int_{\mathbb{R}^3}\f{-e^{-\f {5s}{2}}e^{-s}z_\nu (R_+-R_e)(y-z)}{(e^{-3s}|z_1|^2+e^{-s}|\ck z|^2)^{\f32}}dz.$$
Note that $R_e$ and $R_+$ are supported in $\mathcal X_+$ and $y\in \mathcal X(s)\subset \mathcal X_+$, thus $z\in 2\mathcal X_+$. Hence direct computation shows that
\begin{equation}\label{P1}
\int_{2\mathcal X_+}|P_1(z)|dz=\int_{\substack{|z_1|\leq 2e^{\f 32s}\\|\ck z|\leq2 e^\f s2}}\f{e^{-4s}z_1   }{(e^{-3s}|z_1|^2+e^{-s}|\ck z|^2   )^{\f32} }dz=\int_{|x|\leq2}\f{x_1}{|x|^3}dx\les 1.
\end{equation}
Similarly, we also get
\begin{equation}\label{Pnu}
P_\nu(z)=\f{e^{-3s}z_\nu   }{(e^{-3s}|z_1|^2+e^{-s}|\ck z|^2   )^{\f32} }\in L^1_{2\mathcal{X}_+}.
\end{equation}
Therefore, the derivatives of $\Phi$ can be bounded by
\begin{align*}
&|\p^\ga\p_1\Phi|=\int_{2\mathcal{X}_+}\f{-e^{-\f {5s}{2}}e^{-3s}z_1(\p^\ga R_+-\p^\ga R_e)(y-z)}{(e^{-3s}|z_1|^2+e^{-s}|\ck z|^2)^{\f32}}dz
\les e^{-\f32s}(\|\p^\ga R_e\|_{L^\infty}+\|\p^\ga R_+\|_{L^\infty})
\end{align*}and 
\begin{align*}
&|\p^\ga\p_\nu\Phi|=\int_{2\mathcal{X}_+}\f{-e^{-\f {5s}{2}}e^{-s}z_\nu\p^\ga (\p^\ga R_+-\p^\ga R_e)(y-z)}{(e^{-3s}|z_1|^2+e^{-s}|\ck z|^2)^{\f32}}dz
\les e^{-\f s2}(\|\p^\ga R_e\|_{L^\infty}+\|\p^\ga R_+\|_{L^\infty}).
\end{align*}

\qed

Combining the bootstrap bounds and the above estimates of electric potential, we can obtain the $L^\infty$ bounds of external forces as follows.

\begin{lem}
For $y\in\mathcal{X}(s)$ we have various bounds on the external forces:
\begin{equation}\label{pFWZ}
|\p^\ga F_W|+e^{\f s2}|\p^\ga F_Z|\les
\begin{cases}
e^{-\f s2}, &|\ga|=0\\
e^{-s}\eta^{-\f16+\f{2|\ga|-1}{3(2k-5)}}, &\ga_1\geq1,\ |\ga|=1,2\\
Me^{-s}, & \ga_1=0,\ |\check{\ga}|=1\\
e^{-(1-\f{2}{2k-7})s}, &\ga_1=0,\ |\check{\ga}|=2\\
e^{-(1-\f{3}{2k-7})s}, &|\ga|=3\\
e^{-(1-\f{4}{2k-7})s}, &|\ga|=4
\end{cases}
\end{equation}
and 
\begin{equation} \label{pFUnu}
|\p^\ga F_{U_\nu}|\les
\begin{cases}
e^{-s},\ & |\ga|=0\\
Me^{-s}\eta^{-\f 16},\ & \ga_1=0,\ |\ck\ga|=1\\
e^{-s+\f{3s}{2k-7}}\eta^{-\f16},\ & \ga_1=0,\ |\ck\ga|=2\\
\end{cases}
\end{equation}

\end{lem}
\proof We divide the proof into two steps. 

\textbf{Step 1.} Proof of the bound on $\p^\ga F_W$ and $\p^\ga F_Z$.

Recall (\ref{FW1}),
$$F_W=-2\be_3\be_\tau S(\p_\nu U_\nu)-2\be_1\be_\tau e^s \p_1\Phi.$$
Then by (\ref{pPhi1}), we have
\begin{equation}\label{ppFW}
|\p^\ga F_W|\les |\p^\ga(S\p_\nu U_\nu)|+e^s|\p^\ga\p_1\Phi|\les
|\p^\ga(S\p_\nu U_\nu)|+e^{-\f s2}(\|\p^\ga R_e\|_{L^\infty}+\|\p^\ga R_+\|_{L^\infty}).
\end{equation}

When $|\ga|=0$, then by (\ref{pU1S}), (\ref{pR+}), (\ref{pR}) and (\ref{pUnu}), (\ref{ppFW}) yields
$$|F_W|\les M\ve^\f12e^{-\f s2}+e^{-\f s2}\les e^{-\f s2}.$$

When $|\ga|\leq2$ and $\ga_1\geq1$, using the Leibniz rule, we get
\begin{align*}
|\p^\ga(S\p_\nu U_\nu)|\les
|S\p^\ga(\p_\nu U_\nu)|+|\p^\ga S\p_\nu U_\nu|.
\end{align*}
By (\ref{pU1S}), (\ref{pUnu}) and (\ref{pUnu3}), we have
\begin{align*}
|S\p^\ga(\p_\nu U_\nu)|\les e^{-\f32s}e^{\f{2|\ga|-1}{2k-5}s}\les e^{-s}e^{-\f s2+\f{2|\ga|-1}{2k-5}s}.   
\end{align*}
Then (\ref{eta}) gives
\begin{align*}
|S\p^\ga(\p_\nu U_\nu)|\les \ve^{\f16-\f{2|\ga|-1}{2k-5}}e^{-s} \eta^{-\f16+\f{2|\ga|-1}{3(2k-5)}}\les e^{-s}\eta^{-\f{1}{6}+\f{2|\ga|-1}{3(2k-5)}}.
\end{align*}
Similarly we also have
\begin{align*}
|\p^\ga S\p_\nu U_\nu|\les M^2e^{-s}\ve^{\f12}\eta^{-\f13}
\end{align*}
Meanwhile, (\ref{pR}) and (\ref{pR+}) yield
\begin{align*}
e^{-\f s2}(\|\p^\ga R_e\|_{L^\infty}+\|\p^\ga R_+\|_{L^\infty})\les M^\f1\al e^{-s}\eta^{-\f13}+e^{-2s }.
\end{align*}
Therefore, from (\ref{ppFW}) we obtain $$|\p^\ga F_W|\les e^{-s}\eta^{-\f16+\f{2|\ga|-1}{3(2k-5)}}.$$

Other cases are similar and we just give the sketch of the proof. When $|\ga|=|\ck\ga|=1$,  (\ref{pU1S}), (\ref{pR}) and (\ref{pUnu}) give
$$|\ck\nabla (S\p_\nu U_\nu)|\les |S\ck\nabla^2U_\nu|\les Me^{-s}. $$ (\ref{pR}), and (\ref{pR+}) imply that
$$e^{-\f s2}(\|\p^\ga R_e\|_{L^\infty}+\|\p^\ga R_+\|_{L^\infty})\les e^{-s}.$$
Therefore,  by (\ref{ppFW}) we get
$$|\p_\nu F_W|\les Me^{-s}.$$

When $|\ga|=|\ck\ga|=2$, (\ref{pU1S}),  (\ref{pR}), (\ref{pUnu}) and (\ref{pUnu3}) yield
$$|\ck\nabla^2 (S\p_\nu U_\nu)|\les |S\ck\nabla^3U_\nu|\les e^{-(1-\f{2}{2k-7})s}.$$ Meanwhile, (\ref{pR}) and (\ref{pR+}) imply that
$$e^{-\f s2}(\|\p^\ga R_e\|_{L^\infty}+\|\p^\ga R_+\|_{L^\infty})\les M^\f1\al e^{-s}\eta^{-\f16}.$$
Hence, by (\ref{ppFW}) we get $$|\p^\ga F_W|\les  e^{-(1-\f{2}{2k-7})s}.$$

Finally, When $|\ga|=3,4$, we have 
$$|\p^\ga F_W|\les  e^{-(1-\f{|\ga|}{2k-7})s}.$$

The fact that $F_W=e^{\f s2}F_Z$, combined with the various bounds on $F_W$, gives (\ref{pFWZ}).

\textbf{Step 2.} Proof of the bound on $\p^\ga F_{U_\nu}$.

Recall (\ref{FUnu1}) ,
$$F_{U_\nu}=-2\be_3\be_\tau e^{-\f{s}{2}}S\p_\nu S-2\be_1\be_\tau e^{-\f{s}{2}}\p_\nu\Phi,$$
then by (\ref{pU1S}) and (\ref{pPhick}), we have 
\begin{equation*}
|\p^\ga F_{U_\nu}|\les e^{-\f s2}|\p^\ga(S\p_\nu S)|+e^{-s}(\|\p^\ga R_e\|_{L^
\infty}+\|\p^\ga R_+\|_{L^\infty}).
\end{equation*}

When $|\ga|=0$, by (\ref{pU1S}), (\ref{pR}) and (\ref{pZ}), we have
$$|F_{U_\nu}|\les e^{-s}+ e^{-s}\les e^{-s}.$$

When $\ga_1=0$, $|\ck\ga|=1$, by (\ref{pU1S}), (\ref{pPhick}) and (\ref{pR}),we get
$$e^{-\f s2}|\ck\nabla(S\p_\nu S)|\les e^{-\f s2}|S\ck\nabla^2S|\les Me^{-s}\eta^{-\f16}.$$ Also, (\ref{pR}) and (\ref{pR+}) yield
$$e^{-s}(\|\p^\ga R_e\|_{L^\infty}+\|\p^\ga R_+\|_{L^\infty})\les e^{-\f 32s}.$$
Therefore, we have
$$|\ck\nabla F_{U_\nu}|\les Me^{-s}\eta^{-\f16}.$$

When $\ga_1=0$, $|\ck\ga|=2$,  from (\ref{pU1S}),(\ref{pPhick}) and (\ref{pR}) we have
$$e^{-\f s2}|\ck\nabla^2(S\p_\nu S)|\les e^{-\f s2}|S\ck\nabla^3 S|\les e^{-s+\f{3s}{2k-7}}\eta^{-\f16},$$ 
Meanwhile, $$e^{-s}(\|\p^\ga R_e\|_{L^\infty}+\|\p^\ga R_+\|_{L^\infty})\les M^\f1\al e^{-\f 32s}\eta^{-\f16}.$$
Thus we get
$$|\p^\ga F_{U_\nu}|\les e^{-s+\f{3s}{2k-7}}\eta^{-\f16}.$$

\qed
\begin{cor}
For $y\in\mathcal{X}(s)$ and $k\geq18$, we have 
\begin{equation}\label{pFW2}
|F^{(\ga)}_W|\les
\begin{cases}
e^{-\f s2}, & |\ga|=0\\
\ve^{\f18}\eta^{-\f12+\f{1}{3(2k-5)}}(y), & \ga=(1,0,0)\\
M^\f13\eta^{-\f13}(y), &\ga_1\geq1,\ |\ga|=2\\
\eta^{-\f13}(y), & \ga_1=0,\ |\ck\ga|=1\\
M^\f32\eta^{-\f14}, & \ga_1=0,\ |\ck\ga|=2
\end{cases}
\end{equation}

\begin{equation}\label{pFZ2}
|F^{(\ga)}_Z|\les
\begin{cases}
e^{-s}, & |\ga|=0\\
e^{-\f32s}, &\ga=(1,0,0)\\
e^{-\f32s}, &\ga_1\geq1,\ |\ga|=2\\
Me^{-\f32s}, &\ga_1=0,\ |\ck\ga|=1\\
e^{-\f32s+\f{3s}{2k-7}}, &\ga_1=0,\ |\ck\ga|=2
\end{cases}
\end{equation}

\begin{equation}\label{pFU2}
|F^{(\ga)}_{U_\nu}|\les
\begin{cases}
e^{-s},\ & |\ga|=0\\
M^2e^{-s}\eta^{-\f 16}(y),\ & \ga_1=0,\ |\ck\ga|=1\\
Me^{-s+\f{3s}{2k-7}}\eta^{-\f16}(y),\ & \ga_1=0,\ |\ck\ga|=2\\
\end{cases}
\end{equation}
\end{cor}

\proof
We just need to consider the cases $|\ga|\geq1.$

\textbf{Step 1.} Proof of the bound on $F^{(\ga)}_W$.

Recall (\ref{FW2}),
\begin{align*}
F^{(\ga)}_W=&\p^\ga F_W-\sum_{0\leq\be<\ga}\binom\ga\be\left(\p^{\ga-\be}G_W\p_1\p^\be W+\p^{\ga-\be}h^\mu\p_\mu\p^\be W \right)\notag \\
&-\mathbbm{1}_{|\ga|\geq2}\be_\tau\sum_{\substack{|\be|=|\ga|-1\\ \ga_1=\be_1}}\p^{\ga-\be}W\p^\be \p_1W-\mathbbm{1}_{|\ga|\geq3}\be_\tau\sum_{\be\leq\ga-2}\binom\ga\be \p^{\ga-\be}W\p^\be \p_1W.
\end{align*}

In the following we will utilize the estimates (\ref{pW}), (\ref{pZ}), (\ref{pG}), (\ref{ph}), (\ref{pFWZ}) and (\ref{eta})  repeatedly. When $\ga=(1,0,0)$, we have
\begin{align*}
|F^{(1,0,0)}_W|\les&|\p_1F_W|+|\p_1G_W||\p_1W|+|\p_1h^\mu_W||\p_\mu W|\\
\les&Me^{-s}\eta^{-\f16+\f{1}{3(2k-5)}}+Me^{-s}\eta^{-\f13}+Me^{-2s}\\
\les&M\ve^{\f13}\eta^{-\f12+\f{1}{3(2k-5)}}\\
\les&\ve^\f18\eta^{-\f12+\f{1}{3(2k-5)}}.
\end{align*}

When $\ga=(2,0,0),$ similarly we have 
\begin{align*}
|F^{(2,0,0)}_W|\les &|\p_{11}F_W|+|\p_1G_W||\p_{11}W|+|\p_{11}G_W||\p_1W|+|\p_1h^\mu||\p_{1\mu}W|+|\p_{11}h^\mu||\p_\mu W|+|\p_{11}W||\p_1W|\\
\les&e^{-s}\eta^{-\f16+\f{2|\ga|-1}{3(2k-5)}}+Me^{-s}M^\f13\eta^{-\f13}+Me^{-s}\eta^{-\f13}+Me^{-2s}M^\f23\eta^{-\f13}+e^{-(2-\f{2|\ga|-1}{2k-5})s}+M^\f13\eta^{-\f23}\\
\les &M^\f53e^{-s}\eta^{-\f16+\f{1}{3(2k-5)}}\\
\les & \eta^{-\f13}.
\end{align*}

When $\ga=(1,1,0)$, by  (\ref{pW}), (\ref{pZ}), (\ref{pG}), (\ref{ph}), (\ref{pFWZ}) and (\ref{eta}) again, we get
\begin{align*}
|F^{(1,1,0)}_W|
\les&|\p_{1\nu}F_W|+|\p_1G_W||\p_{1\nu}W|+|\p_\nu G_W||\p_{11}W|+|\p_{1\nu}G_W||\p_1W|\\
&+|\p_1h^\mu||\p_{1\mu}W|+|\p_\nu h^\mu||\p_{1\mu}W|+|\p_{1\nu}h^\mu||\p_\mu W|+|\p_\nu W||\p_{11}W|\\
\les & e^{-s}\eta^{-\f16+\f{2|\ga|-1}{3(2k-5)}}+Me^{-s}M^{\f23}\eta^{-\f13}+M\ve^{\f12}M^\f13\eta^{-\f13}+Me^{-s}\eta^{-\f13}\\
&+Me^{-2s}M^\f23\eta^{-\f13}+M\ve^{\f12}e^{-s}M^\f23\eta^{-\f13}+Me^{-\f32s}+M^\f13\eta^{-\f13}\\
\les &M^\f13\eta^{-\f13}.
\end{align*}

When $\ga=(0,1,0)$, we also have
\begin{align*}
	|F^{(0,1,0)}_W|\les& |\ck\nabla F_W|+|\ck\nabla G_W||\p_1W|+|\ck\nabla h^\mu||\p_\mu W|\\
	\les&Me^{-s}+M\ve^\f12\eta^{-\f13}+M\ve^\f12e^{-s}\\
	\les&\eta^{-\f 13}.
\end{align*}

When $\ga=(0,2,0)$ and $k\geq18$, it holds that
\begin{align*}
	|F^{(0,2,0)}|\les&|\ck\nabla^2F_W|+|\ck\nabla^2G_W||\p_1W|+|\ck\nabla  G_W||\p_1\ck\nabla W|+|\ck\nabla^2h^\mu||\p_\mu W|+|\ck\nabla h^\mu|| \ck\nabla^2W|+|\ck\nabla W||\p_1\ck\nabla W|\\
	\les& e^{-(1-\f{2}{2k-7})s}+Me^{-\f s2}\eta^{-\f13}+M^\f53\ve^\f12\eta^{-\f13}+Me^{-\f32s}+M\ve^\f12e^{-s}+M^2\ve^\f12e^{-s}\eta^{-\f16}+M^\f23\eta^{-\f13}\\
	\les& M^\f23\eta^{-\f14}.
\end{align*}
	
\textbf{Step 2.} Proof of the bounds on $F^{(\ga)}_Z$.

Recall from (\ref{FZ2}) that
\begin{align*}
	F^{(\ga)}_Z=&\p^\ga F_Z-\sum_{0\leq\be<\ga}\binom\ga\be\left(\p^{\ga-\be}G_Z\p_1\p^\be Z+\p^{\ga-\be}h^\mu\p_\mu\p^\be Z \right)\notag\\
	&-\be_2\be_\tau\sum_{\substack{|\be|=|\ga|-1\\ \ga_1=\be_1}}\p^{\ga-\be}W\p^\be \p_1Z-\mathbbm{1}_{|\ga|\geq2}\be_2\be_\tau\sum_{\be\leq\ga-2}\binom\ga\be \p^{\ga-\be}W\p^\be \p_1Z.
\end{align*}
 We  keep utilizing the estimates (\ref{pW}), (\ref{pZ}), (\ref{pG}), (\ref{ph}) and (\ref{pFWZ}), repeatedly. 

When $\ga=(1,0,0),$ we have 
\begin{align*}
|F^{(1,0,0)}_Z|\les& |\p_1F_Z|+|\p_1G_Z||\p_1Z|+|\p_1h^\mu_Z||\p_\mu Z|+|\p_1W||\p_1Z|\\
\les & e^{-\f32s}\eta^{-\f16+\f{1}{3(2k-5)}}+M^2e^{-s}e^{-\f32s}+M^2e^{-2s}\ve^\f12e^{-\f s2}+M\eta^{-\f13}e^{-\f32s}\\
\les &e^{-\f32s}.
\end{align*}

When $\ga_1\geq1,$ $|\ga|=2$, it holds that
\begin{align*}
	|F^{(\ga)}_Z|\les&|\p_1\nabla F_Z|+|\p_1\nabla G_Z||\p_1Z|+|\p_1G_Z||\p_1\nabla Z|+|\nabla G_Z||\p_{11}Z|\\
	&+|\p_1\nabla h^\mu||\p_\mu Z|+|\p_1h^\mu||\p_\mu\nabla Z|+|\nabla h^\mu||\p_1\p_\mu Z|\\
	&+|\ck\nabla W||\p_{11}Z|+|\p_1W||\p_{11}Z|+|\p_1\nabla W||\p_1Z|\\
	\les&e^{-\f32s}\eta^{-\f16+\f{1}{3(2k-5)}}+M^2e^{-\f52s}+M^2e^{-\f52s}+M^2e^{-\f52s}\\
	&+M^2\ve^\f12e^{-\f52s}+M^2e^{-3s}+M^2\ve^\f12e^{-\f52s}\\
	&+Me^{-\f32s}+Me^{-\f32s}\eta^{-\f13}+M^{\f53}e^{-\f32s}\eta^{-\f13}\\
	\les&e^{-\f32s}.
\end{align*}

When $\ga=(0,1,0)$, we have
\begin{align*}
|F^{(0,1,0)}_Z|\les&|\ck\nabla F_Z|+|\ck\nabla G_Z||\p_1Z|+|\ck h^\mu||\p_\mu Z|+|\ck\nabla W||\p_1Z|\\
\les&Me^{-\f 32s}+M^2\ve^\f12e^{-\f32s}+M^2\ve e^{-\f32s}+Me^{-\f32s}\\
\les&Me^{-\f 32s}.
\end{align*}

When $\ga_1=0,$ $|\ck\ga|=2$, we get
\begin{align*}
|F^{(0,2,0)}_Z|\les&|\ck\nabla^2F_Z|+|\ck\nabla^2G_Z||\p_1Z|+|\ck\nabla G_Z||\p_1\ck\nabla Z|+|\ck\nabla^2h^\mu||\p_\mu Z|+|\ck\nabla^2h^\mu||\ck\nabla^2Z|\\
&+|\ck\nabla W||\p_1\ck\nabla Z|+|\ck\nabla^2W||\p_1Z|\\
\les&e^{-(\f32-\f{2}{2k-7})s}+M^2e^{-3s}+M^2\ve^\f12e^{-\f s2}+M^2e^{-\f52s}+Me^{-\f32s}+M^2e^{-\f32s}\eta^{-\f16}\\
\les&e^{-(\f32-\f{2}{2k-7})s}.
\end{align*}

\textbf{Step 3.} Proof of the bounds on $F^{(\ga)}_{U_\nu}$.

The proof of bounds for  $F^{(\ga)}_{U_\nu}$ are very similar to those of $F^{(\ga)}_W$ and $F^{(\ga)}_Z$, and we left it to readers.

\qed
\begin{lem}
For $\p^\ga \td{F}_W$ and $|y|\leq L$, we have 
\begin{equation}\label{ptdFW}
|\p^\ga \td{F}_W|\les \ve^\f{1}{10}
\begin{cases}
\eta^{-\f16}(y), & |\ga|=0\\
e^{-s}\eta^{-\f16+\f{2|\ga|-1}{3(2k-5)}}(y), & \ga=(1,0,0)\\
\eta^{-\f16}(y), & \ga_1=0,\ \  |\ck\ga|=1\\
1 & |\ga|=4,\ \ |y|\leq l
\end{cases}
\end{equation}
Furthermore, for $y=0$ and $|\ga|=3$ we have
\begin{equation}\label{p3tdFW}
|(\p^\ga\td{F}_W)^0|\les e^{-(\f12-\f{4}{2k-7})s}.
\end{equation}
\end{lem}
\proof
Due to (\ref{G}), we have
 \begin{equation}\label{Gve}
|G_W|\leq M(e^{-\f s2}+Le^{-s}+L\ve^\f12)\leq\ve^{\f18}.
\end{equation}for $|y|\leq L.$
Recall from \eqref{FW3}, we have $$\td{F}_W=F_W-e^{-\f s2}\be_\tau \dot{\ka}+((\be_\tau-1)\bar{W}-G_W)\p_1\bar{W}-h^\mu\p_\mu\bar{W}.$$

When $|\ga|=0$, by (\ref{betau}) and (\ref{pbarW}) we get
\begin{align*}
|\td{F}_W|\les&|F_W|+Me^{-\f s2}+|((\be_\tau-1)\bar{W}+G_W)||\p_1\bar{W}|+|h^{\mu}||\p_\mu W|\\
\les & e^{-\f s2}+Me^{-\f s2}+(M\ve\eta^\f16+\ve^\f18)\eta^{-\f13}+M\ve^\f12e^{-\f s2}\\
\les & \ve^\f{1}{10}\eta^{-\f16}.
\end{align*}

When $\ga=(1,0,0)$, (\ref{betau}) gives
\begin{align*}
|\p_1\td{F}_W|\les&|\p_1F_W|+(\be_\tau-1)|\p_1\bar{W}|^2+|\p_1G_W||\p_1\bar{W}|+|(\be_\tau-1)\bar{W}+G_W||\p_{11}\bar{W}|+|\p_1h^\mu||\p_\mu \bar{W}|+|h^\mu||\p_{1\mu}\bar{W}|\\
\les& e^{-s}\eta^{\f16+\f{2|\ga|-1}{3(2k-5)}}+M\ve\eta^{-\f23}+Me^{-\f s2}\eta^{-\f13}+(M\ve\eta^\f16+\ve^\f18)\eta^{-\f23}+Me^{-2s}+M\ve^\f12e^{-\f s2}\eta^{-\f13}\\
\les &\ve^\f{1}{10} \eta^{-\f12+\f{1}{3(2k-5)}}.
\end{align*}

When $\ga=(0,1,0)$, it holds from (\ref{pbarW}) and (\ref{pG}) that
\begin{align*}
|\p_\nu\td{F}_W|\les&|\p_\nu F_W|+(\be_\tau-1)|\p_\nu\bar{W}||\p_1\bar{W}|+|\p_\nu G_W||\p_1\bar{W}|+|(\be_\tau-1)\bar{W}+G_W||\p_{1\nu}\bar{W}|+|\p_\nu h^\mu||\p_\mu \bar{W}|+|h^\mu||\p_{\nu\mu}\bar{W}|\\
\les& Me^{-s}+M\ve\eta^{-\f13}+M\ve^\f12\eta^{-\f13}+(M\ve\eta^\f16+\ve^\f18)\eta^{-\f13}+Me^{-\f s2}\eta^{-\f13}+M\ve^\f12e^{-\f s2}\eta^{-\f16}\\
\les& \ve^\f{1}{10}\eta^{-\f16}.
\end{align*}

When $|\ga|=4$, it holds for $|\be|\leq4$ that $$|\p^\be\p_1(\bar W)^2|\les1,\ \ \ |\p^\be\p_1\bar W|\les \eta^{-\f13},\ \ \ |\p^\be\p_\mu\bar W|\les 1.$$
Hence, from (\ref{pbarW}), (\ref{pG}) and (\ref{ph}) we get
\begin{align*}
|\nabla^4\td F_W|\les&|\nabla^4 F_W|+|\nabla^4\Big(((\be_\tau-1)\bar W-G_W)\p_1\bar W\Big)|+|\nabla^4(h^\mu\p_\mu\bar W)|\\
\les&e^{-(1-\f{4}{2k-7})s}+M\ve|\nabla^4(\p_1\bar W)^2|+\sum_{\be\leq\ga}|\p^{\ga-\be} G_W||\p^\be\p_1\bar W|+\sum_{\be\leq\ga}|\p^{\ga-\be} h^\mu||\p^\be\p_\mu\bar W|\\
\les& \ve^\f14+(\ve^\f18+M\ve^\f12)\eta^{-\f13}+M\ve^\f12e^{-\f s2}\\
\les& \ve^\f{1}{10}.
\end{align*}

When $|\ga|=3$, we note that the even order derivatives of $\bar{W}$ and $\ck\nabla\bar{W}$ vanish at $y=0$, as well as $\p_1\bar W(0)=-1$. Therefore, from the expression of $\td F_W$, according to Leibniz rule , combined with (\ref{pFWZ}), (\ref{pbarW}), (\ref{pG}) and (\ref{ph}) we have
\begin{align*}
|(\nabla^3\td{F}_W)^0|\les&|(\nabla^3F_W)^0|+|(\nabla^3((\be_\tau-1)\bar{W}-G_W))^0|+|(\nabla((\be_\tau-1)\bar{W}-G_W))^0|+|(\nabla h^\mu)^0|\\
\les&|(\nabla^3F_W)^0|+|1-\be_\tau|+|(\nabla^3G_W)^0|+|(\nabla G_W)^0|+|(\nabla h^\mu)^0|\\
\les&e^{-(1-\f{3}{2k-7})s}+Me^{-s}+e^{-(1-\f{2}{2k-7})s}+Me^{-s}+e^{-s}\\
\les&M\ve^{\f{1}{2k-7}}e^{-(1-\f{3}{2k-7})s}\\
\les&\ve^{\f{1}{10}}.
\end{align*}

\qed

\begin{cor}
For $\td{F}^{(\ga)}_W$ and $|y|\leq L$ we have
\begin{equation}\label{ptdFW2}
|\td{F}^{(\ga)}_W|\les \ve^\f{1}{15}
\begin{cases}
\eta^{-\f12+\f{1}{3(2k-5)}}(y) & \ga=(1,0,0), |y|\leq L\\
\eta^{-\f16}(y) &\ga_1=0, |\ck\ga|=1, |y|\leq L \\

1 &|\ga|=4, |y|\leq l
\end{cases}
\end{equation}
Furthermore, we have
\begin{equation} \label{ptdFW3}
|\td{F}_W^{(\ga),0}|\les e^{-(\f12-\f{4}{2k-7})s}\ \ \ |\ga|=3
\end{equation}
\end{cor}

The proof is  similar to that of Lemma 5.5 and we omit it.

\section{Proof of Proposition 4.1}

\subsection{Velocity Fields in Terms of Lagrangian Trajectories}
We define $\Psi^{y_0}_W$ to be the Lagrangian flow associated with the velocity field $\mathcal V_W$:
\begin{equation*}
\begin{cases}
\p_s \Psi^{y_0}_W(s)=\mathcal V_W(\Psi^{y_0}_W(y,s),s),\\
\Psi^{y_0}(s_0)=y_0.
\end{cases}
\end{equation*}
 We also define $\Psi^{y_0}_Z$ and $\Psi^{y_0}_U$ similarly. In order to close the bootstrap argument, we need to control the integration $\int\eta^{-\de}\circ\Psi$. As in \cite{Buckmaster2019}, it's necessary to analyze the lowest speed of velocity fields $(\mathcal{V}_W, \mathcal{V}_Z,\mathcal{V}_U)$ escaping from the origin.

\begin{lem}
Let $|y_0|\geq l$, $s\geq-\log\ve$, then the Lagrangian flow $\Psi^{y_0}_W$ escapes from the origin at an exponential speed, that is 
\begin{equation}\label{Psispd}
|\Psi^{y_0}_W(s)|\geq|y_0|e^{\f{s-s_0}{5}}.
\end{equation}
Moreover,  for fixed $\de>0$ and $M$ large enough it holds that
\begin{subequations}
\begin{align}
\int_{s_0}^{s}(1+\Psi_W^{y_0}(s'))^{-\de}ds'\leq& 6\log\f1l\ \ \ \ \  &l\leq|y_0|\leq L\label{int1} \\
\int_{s_0}^{s}(1+\Psi_W^{y_0}(s'))^{-\de}ds'\leq& \f{5}{2\de}L^{-2\de}=\f{5}{2\de}\ve^\f\de5  &|y_0|\geq L \label{int2}
\end{align}
\end{subequations}
\end{lem}
\proof
The proof of (\ref{Psispd}) is similar to that in \cite{Buckmaster2019}. It suffices to show
\begin{equation}\label{63}
	y\cdot\mathcal V_W(y)\geq\f15|y|^2,\ \ \ \ |y|\geq l,
\end{equation}which follows from (\ref{pW}), (\ref{ptdW}) and (\ref{G}). 

For (\ref{int1}), since $2\eta(y)\geq(1+|y|^2)$, \eqref{63} implies that $le^{\f15(s'-s_0)}=p$,
\begin{align*}
\int_{s_0}^{s}(1+\Psi_W^{y_0}(s'))^{-\de}ds'\leq&\int_{s_0}^{\infty}(1+l^2e^{\f25(s'-s_0)})^{-\de}ds'.
\end{align*}
Write $le^{\f15(s'-s_0)}=p$, then we have
\begin{align*}
 \int_{s_0}^{s}(1+\Psi_W^{y_0}(s'))^{-\de}ds'\leq&\int_l^\infty\f{5dp}{(1+p^2)^{\de}p}\\
=&\int_l^1\f{5dp}{(1+p^2)^{\de}p}+\int_1^\infty\f{5dp}{(1+p^2)^{\de}p}\\
\leq&5\log\f1l +\f{5}{2\de}\leq6\log\f1l,\\
\end{align*}where we take $M$ large enough in the last inequality and note that $l=M^{-\f{1}{500}}$.

The bound (\ref{int2}) can be established similarly. In fact, from \eqref{Psispd} and $|y_0|\geq L$, we get $$\int_{s_0}^{s}(1+\Psi_W^{y_0}(s'))^{-\de}ds'\leq\int_{s_0}^{\infty}(1+L^2e^{\f25(s'-s_0)})^{-\de}ds'\leq \f{5}{2\de}L^{-2\de}\\
= \f{5}{2\de}\ve^{\f\de5}.$$

\qed

\begin{lem}Let $0\leq\sigma_1<\f12$, $2\sigma_1<\sigma_2$. Let $\Psi(s)=(\Psi_1(s), \ck\Psi(s))$ denote either $\Psi^{y_0}_Z(s)$ or $\Psi^{y_0}_U(s)$. Suppose 
\begin{equation}\label{ka0}
\ka_0\geq\f{3}{1-\max\{\be_1,\be_2\}},
\end{equation}
then for any $y_0\in\mathcal X_0$
\begin{equation}\label{etaint}
\int_{-\log\ve}^{\infty}e^{\sigma_1s'}(1+|\Psi_1(s')|)^{-\sigma_2}ds'\leq C_{\sigma_1,\sigma_2}.
\end{equation}
In particular, 
\begin{align}\label{p1Wint}
\sup_{y_0\in\mathcal{X}_0}\int_{-\log\ve}^{s}|\p_1W|\circ\Psi^{y_0}(s')ds'&\les1.
\end{align}
\end{lem}

\proof
We want to show if $\Psi(s)=\Psi_Z^{y_0}$ or $\Psi^{y_0}_U(s)$, then 
\begin{equation}\label{psi1}
	\f{d}{ds}\Psi_1(s)\leq-\f12e^{\f s2}\ \   \text{if}\ \ \Psi_1\leq e^{\f s2},\ \ s\geq-\log\ve.
\end{equation}
We just verify \eqref{psi1} for  $\Psi(s)=\Psi_Z^{y_0}$ and the other case is similar. By definition, it holds that $$\f{d}{ds}\Psi_1=\f32\Psi_1+\be_2\be_\tau W\circ\Psi+G_U\circ\Psi.$$ Since $\be_2\leq1$, by (\ref{pW}) and taking $\ve$ sufficiently small, under the assumption $\Psi_1(s)\leq e^{\f s2}$, we have
\begin{align*}
	\f{d}{ds}\Psi_1\leq &\f32e^{\f s2}+2\eta^{\f16}(\Psi)-(1-\be_2)\ka_0e^{\f s2}+\ve^\f12e^{\f s2}\\
	\leq &\f32 e^{\f s2}-(1-\be_2)\ka_0 e^{\f s2}+\ve^\f18e^{\f s2}\\
	\leq &-\f12e^{\f s2},
\end{align*} 
where the last inequality we have used $1-\be_2>0$ and (\ref{ka0}).

Next, the following proof is similar to \cite{Buckmaster2019}, which we just sketch it.

Now we claim that 
\begin{equation}\label{e2s}
	|\Psi^{y_0}_1(s)|\geq \min\{ |e^{\f s2}-e^{\f {s_*}{2}}|, e^{\f s2} \}
\end{equation}
for all $y_{1_0}$ starting from $|y_1|\leq\ve^{-1}.$

Then to prove 
 (\ref{etaint}), we first note that since $\int_{-\log\ve}^\infty e^{(\si_1-\f{\si_2}{2})s'}ds'\les1$,
  it suffice to prove 
 $$\int_{-\log\ve}^{\infty}e^{\si_1s'}\Big( 1+|e^{\f{s'}{2}}-e^{\f{s_*}{2}}|  \Big)^{-\si_2}ds'\leq C,$$
which can be obtained by change of variable $r=e^{\f {s'}{2}}$.

\qed

\subsection{$L^\infty$ Estimates for Specific Vorticity $\Omega$}
\begin{lem}
The sound speed $S$ is bounded blow from $0$ strictly, namely 
\begin{equation}\label{Sbd}
\|S(\cdot,s)-\f{\ka_0}{2}\|_{L^\infty_{\mathcal X(s)}}\leq \ve^\f18
\end{equation}for all $s\geq-\log\ve$.
\end{lem}
\proof
Since $S=\f12(e^{-\f s2}W+\ka-Z)$ and $y\in\mathcal X(s)$, (\ref{pW}) and (\ref{pZ}) yield
\begin{align*}
|S-\f{\ka_0}{2}|\leq&|\f{\ka-\ka_0}{2}|+\f{e^{-\f s2}}{2}|W|+\f{|Z|}{2}\\
\leq& \f12\int_{-\ve}^{T_*}|\dot\ka|+2\ve^\f16\eta^{-\f16}\eta^\f16+\f12M\ve^\f12\\
\leq& \ve^\f18,
\end{align*}where in the last inequality we have utilized (\ref{modu}), (\ref{eta}) and $T\leq 3M\ve^2$.

\qed

\begin{lem}
For specific vorticity $\td\zeta(x,t)=\Omega(y,s)$, we have 
\begin{equation}\label{tdze}
\|\td\zeta(\cdot,t)\|_{L^\infty}=\|\Omega(\cdot,s)\|_{L^\infty}\leq2
\end{equation}
\end{lem}
\proof
From (\ref{zeta}), $x=\xcal-\xi(t)$, $t=\f{1+\al}{2}\tcal$, $\td\zeta(x,t)=\zeta(\xcal,t)$ we have 
\begin{equation}\label{evoze}
\f{1+\al}{2}\p_t\td\zeta+((\td u-\dot\xi)\cdot\nabla)\td\zeta-(\td\zeta\cdot\nabla)\td u=0.
\end{equation}

We estimate $\td\zeta_1$ first. (\ref{pUnu}) gives
$$|\td n_e\td\zeta_1|=|(\nabla\times\td u)_1|=|\p_{x_2}\td u_3-\p_{x_3}\td u_2|\leq M\ve^\f12.$$
Meanwhile, $\td n_{e}=(\al\td\sigma)^{\f1\al}$ and (\ref{Sbd}), we have 
$|n_e|\geq \ka_0^\f1\al$. Therefore, we have
\begin{equation}\label{ze1}
|\td{\zeta}_1|\leq \ve^\f13.
\end{equation}

In order to estimate $\td{\zeta}_\nu$, we rewrite the equation (\ref{evoze}),
\begin{equation}\label{tdzenu}
\p_t\td\zeta_\nu+2\be_1((\td u-\dot\xi)\cdot\nabla)\td\zeta_\nu=2\be_1(\td\zeta_1\p_{x_1}\td u_\nu+\td\zeta_2\p_{x_2}\td u_\nu+\td\zeta_3\p_{x_3}\td u_\nu).
\end{equation}
Now for $x_0\in\XX$ we define the trajectory $\psi^{x_0}$ by $\p_t\psi^{x_0}=2\be_1 (\td u-\dot\xi)\circ\psi^{x_0}$ for $t>-\ve$ and $\psi^{x_0}(-\ve)=x_0$. Let $\mathcal Q_j=\td\zeta_j\circ \psi^{x_0}$, then (\ref{tdzenu}) becomes 
$$\p_t\mathcal Q_\nu=2\be_1(\p_{x_\nu}\td u_j\circ\psi^{x_0})\mathcal Q_j.$$ 
Since (\ref{pUnu}) and (\ref{ze1}) lead to $|2\be_1(\p_1\td u_\nu\circ\psi^{x_0})\mathcal Q_1|+|2\be_1\p_{x_\nu}\td u_\mu\circ\psi^{x_0}|\leq M\ve^\f13,$ we have
$$\f12\f{d}{dt}(\mathcal Q^2_2+\mathcal Q^2_3)\leq M\ve^\f13(\mathcal Q^2_2+\mathcal Q^2_3)+M\ve^\f13.$$
Integration from $-\ve$ to $t\in(-\ve, T_*)$, by (\ref{iniOme}) and Gr\"{o}nwall inequality, we have 
$$\mathcal Q^2_2+\mathcal Q^2_3\leq \f32,$$ and (\ref{tdze}) follows.

\qed

\subsection{Closure of $L^\infty$ Type Bootstrap on $Z$}

\begin{lem}Let $\Psi=\Psi_Z^{y_0}$, then we have
\begin{subequations}\label{Zest}
\begin{align}
|Z\circ\Psi(s)|&\les \ve^\f12, \\
e^{\f32s}|\p^\ga Z\circ\Psi(s)|&\les1,\ \ \ga_1\geq1,|\ga|=1,2\\
e^{\f s2}|\ck\nabla Z\circ\Psi(s)|&\les\ve^{\f12},\\
e^{s}|\ck\nabla^2Z\circ\Psi(s)|&\les1.
\end{align}
\end{subequations}
\end{lem}
\proof
For $\mu\leq\f{3\ga_1+\ga_2+\ga_3}{2}$, multiplying $e^{\mu s}$ on both sides of (\ref{Devo2}), we get
\begin{align*}
\p_s(e^{\mu s}\p^\ga Z)+D^{(\ga,\mu)}_Z(e^{\mu s}\p^\ga Z)+(\mathcal{V}_Z\cdot\nabla)(e^{\mu s}\p^\ga Z)=e^{\mu s}F_Z^{(\ga)},
\end{align*}
where $$D_Z^{(\ga,\mu)}=-\mu+\f{3\ga_1+\ga_2+\ga_3}{2}+\be_2\be_\tau\ga_1\p_1W.$$
Composing the Lagrangian trajectories $\Psi_Z$, by Gr\"onwall inequality, we have
\begin{align}
e^{\mu s}|\p^\ga Z\circ\Psi^{y_0}_Z(s) |\leq& \ve^{-\mu}|\p^\ga Z(y_0,-\log\ve)|\exp\left(-\int_{-\log\ve}^{s}D^{(\ga,\mu)}_Z\circ \Psi^{y_0}_Z(s')ds' \right)\notag\\
&+\int_{-\log\ve}^{s}e^{\mu s'}|F^{(\ga)}_Z\circ\Psi^{y_0}_Z(s')|\exp\left(-\int_{s'}^{s}D^{(\ga,\mu)}_Z\circ \Psi^{y_0}_Z(s'')ds'' \right)\label{gronZ}
\end{align}
Note that 
\begin{align*}
&-\int_{-\log\ve}^{s}D^{(\ga,\mu)}_Z\circ \Psi^{y_0}_Z(s')ds'\\
\les&\int_{-\log\ve}^{s}-\f{3\ga_1+\ga_2+\ga_3}{2}+\mu+\int_{-\log\ve}^{s}|\p_1W|ds,
\end{align*}
so by (\ref{p1Wint}) and $\mu\leq\f{3\ga_1+\ga_2+\ga_3}{2}$, we have
\begin{equation}
\exp\left(-\int_{-\log\ve}^{s}D^{(\ga,\mu)}_Z\circ \Psi^{y_0}_Z(s')ds' \right)\les \exp\left(\int_{-\log\ve}^{s}-\f{3\ga_1+\ga_2+\ga_3}{2}+\mu\right)\les1. \label{expZ}
\end{equation}
Combining (\ref{gronZ}) and (\ref{expZ}), we get
\begin{equation}\label{616}
e^{\mu s}|\p^\ga Z\circ\Psi^{y_0}_Z(s) |\leq \ve^{-\mu}|\p^\ga Z(y_0,-\log\ve)|
+\int_{-\log\ve}^{s}e^{\mu s'}|F^{(\ga)}_Z\circ\Psi^{y_0}_Z(s')|.
\end{equation}

Next we plug initial data (\ref{pZ}) and external forces bounds (\ref{pFZ2}) into (\ref{616}) for different $\ga$.
When $|\ga|=0$, we let $\mu=0,$ and then (\ref{616}) becomes
\begin{align*}
| Z\circ\Psi^{y_0}_Z(s)| \les \ve+\int_{-\log\ve}^{s}e^{-s'}ds'\les 2\ve\les\ve^\f12.
\end{align*}
When $\ga=(1,0,0)$, we let $\mu=\f32$. Then \eqref{616} yields
\begin{align*}
|\p_1Z\circ\Psi^{y_0}_Z(s)| \les& \ve^{-\f32 s}|\p_1Z(y_0,-\log\ve,)|+\int_{-\log\ve}^{s}e^{\f32s'}|F^{(1,0,0)}_Z\circ\Psi^{y_0}_Z(s')|ds'\\
\les&1+\int_{-\log\ve}^{s}\left(1+|\Psi_1(s')|^2 \right)^{-\f{1}{3(2k-5)}}ds'\\
\les&1,
\end{align*}
where in the last inequality we have used (\ref{etaint}) for $\sigma_1=0$, $\sigma_2=\f{1}{3(2k-5)}$.

When $\ga_1\geq1,$ $|\ga|=2$, we let $\mu=\f32$. Then $$\mu-\f{3\ga_1+\ga_2+\ga_3}{2}\leq-\f12.$$ So we have
\begin{align*}
e^{\f32s}|\p^\ga Z\circ\Psi_Z^{y_0}|\les&\ve^{-\f32}|\p^\ga Z(y_0, -\log\be)|+\int_{-\log\ve}^s e^{\f32 s'} e^{-\f32 s'}e^{-\f12(s-s')}ds'\les 1.
\end{align*}

At last when $|\ga|=|\ck\ga|=1$ or $2$, after setting $\mu=\f12$ and $1$, we have 
\begin{align*}
e^{\f s2}|\ck\nabla Z\circ\Psi|\les& \ve^{-\f12}|\ck\nabla Z(y_0,-\log\ve)|+M^2\int_{-\log\ve}^{s}e^{-s'}ds' \les\ve^{\f12}+M^2\ve\les\ve^{\f12}
\end{align*}
and
\begin{align*}e^{ s}|\ck\nabla^2 Z\circ\Psi|\les& \ve^{-\f12}|\ck\nabla^2 Z(y_0,-\log\ve)|+M\int_{-\log\ve}^{s} e^{-\f{s'}{2}+\f{3s'}{2k-7}}ds'\les 1,
\end{align*}
for $k\geq 18.$
\qed

\subsection{Closure of $L^\infty$ Type Bootstrap on $U_\nu$}

\begin{lem}
For $\Psi=\Psi^{y_0}_U$, we have
\begin{subequations}\label{Uest}
\begin{align}
|U_\nu\circ\Psi(s)|&\les \ve^{\f12},\\
e^{\f s2}|\ck\nabla U_\nu\circ\Psi(s)|&\les \ve^{\f12},\\
e^s |\ck \nabla^2U_\nu\circ\Psi(s)|&\les1.
\end{align}
\end{subequations}
Besides, it holds that\begin{equation}
e^{\f32s}|\p_1 U_\nu|\les1.
\end{equation}
\end{lem}
\proof
The proof of (\ref{Uest}) is very similar to Lemma 6.4 by using (\ref{iniUnu}), (\ref{pUnu}) and (\ref{pFUnu}) so we omit it.

 For estimating $\p_1 U_\nu$, we use the specific vorticity, $\nabla\times u=\rho \zeta=(\al\sigma)^{\f1\al}\zeta$. Rewrite it to self-similar variables we get
$$|e^{\f32s}\p_1 U_2-e^{\f s2}\p_2 U_1|\les| (\al S)^{\f1\al}\Omega|.$$ Then by (\ref{pU1S}) and (\ref{tdze}) we have
$$e^{\f32 s}|\p_1 U_2|\les \ka_0^\f1\al+1\les1.$$
 And it is similar to establish $|\p_1 U_3|\les e^{-\f32s}.$

\qed

\subsection{Closure of $L^\infty$ Type Bootstrap on $\td W$}

\begin{lem}
We have
\begin{subequations}
\begin{align}
|\p^\ga \td{W}|\leq &\f12M\ve^{\f{1}{18}}l^{4-|\ga|}                      \ &|\ga|\leq3,\ |y|\leq l \\
|\p^\ga \td{W}|\leq& \f12M\ve^{\f{1}{18}} &|\ga|=4,\ |y|\leq l
\end{align}
\end{subequations}
\end{lem}

\proof
\textbf{Case 1.} $|\ga|=4,\ |y|\leq l$

First recall \ref{DtdW}
 \begin{align}
	\left(\p_s+\f{3\ga_1+\ga_2+\ga_3-1}{2}+\be_\tau(\p_1\bar{W}+\ga_	1\p_1W ) \right)\p^\ga\td{W}+(\mathcal{V}_W\cdot\nabla)\p^\ga\td{W}=\td{F}^{(\ga)}_W,
\end{align}
In this case, the damping term are bounded from below
\begin{align*}
D^{(\ga)}_{\td{W}}=&\f{3\ga_1+\ga_2+\ga_3-1}{2}+\be_\tau(\p_1\bar{W}+\ga_1\p_1W)\\
=&\f32+\ga_1+\be_\tau(\p_1\bar{W}+\ga_1\p_1W)\\
\geq&\f32+\ga_1-(1+2M^2\ve)(1+\ga_1(1+\ve^{\f{1}{20}}))\\
\geq&\f13,
\end{align*}
where we have used (\ref{pbarW}), (\ref{betau}) and (\ref{pWlbd}) in the first inequality.

Note that if $\Psi^{y_0}_W(s)$ denotes the flow generated by velocity field $\mathcal{V}_W$, then we have
$$
\f{d}{ds}\left(\p^\ga\td{W}\circ\Psi^{y_0}_W \right)+\left(D^{(\ga)}_{\td{W}}\circ\Psi^{y_0}_W\right) \left(\p^\ga\td{W}\circ\Psi^{y_0}_W \right)=\td{F}^{(\ga)}_{\td{W}}\circ\Psi^{y_0}_W.
$$
Then applying Gr\"onwall's inequality as well as (\ref{pFW2}), we  get 
\begin{equation}
|\p^\ga\td{W}\circ\Psi^{y_0}_W|\les \ve^{\f{1}{15}} +|\p^\ga\td{W}(y_0,-\log\ve)|\les \ve^{\f{1}{15}}+\ve^{\f14}\leq \f12M \ve^{\f{1}{18}}.
\end{equation}

\textbf{Case 2.} $|\ga|\leq3$, $|y|\leq l$.

By (\ref{conW}), $\td{W}=W-\bar{W}$ and the explicit expression of $\bar{W}$, we have
\begin{equation}\label{detdW}
\td{W}(0,s)=\nabla\td{W}(0,s)=\nabla^2\td{W}(0,s)=0.
\end{equation}
By (\ref{DtdW}), when $y=0$ there holds 
$$
\p_s(\p^\ga\td{W})^0=\td{F}_W^{(\ga),0}-G^0_W(\p_1\p^\ga\td{W})^0-h^{\mu,0}(\p_\mu\p^\ga\td{W})^0-(1+\ga_1)(1-\be_\tau)(\p^\ga\td{W})^0.
$$
Then by (\ref{ptdFW3}), (\ref{ptdW4}) and (\ref{Gh0}), we have
$$|\p_s(\p^\ga\td{W})^0|\les e^{-(\f12-\f{4}{2k-7})s}+Me^{-s}\ve^{\f{1}{18}}+Me^{-\f s2}\ve^{\f{1}{18}}+Me^{-s}\ve^{\f{1}{18}}\les  e^{-(\f12-\f{4}{2k-7})s}. $$
Thus, it follows from Newton-Leibniz formula that
$$|\p^\ga\td{W}(0,s)|\leq |\p^\ga\td{W}(0,-\log\ve)|+\int_{-\log\ve}^{s}|\p_s(\p^\ga W)^0(s')|ds'\leq\f{1}{10}\ve^{\f{1}{18}}.$$ 
Thus, we get from \eqref{detdW} that
\begin{equation}\label{detdW2}
|\p^\ga\td{W}(0,s)|\leq\f{1}{10}\ve^{\f{1}{10}},\ \ \ |\ga|\leq3,
\end{equation}

At last, 
when $|y|\leq l,$ according to Taylor expansion, it is easy to get $|\p^\ga\td W(y,s)|\leq \f12 M\ve^\f{1}{18}l^{4-|\ga|}$ by integrating from $0$
to $y.$

\qed

\begin{lem}
For $|y|\leq L$, we have 
\begin{subequations}
\begin{align*}
|\td{W}\circ\Psi_W^{y_0}(s)|\leq \f12\ve^\f{1}{20} \eta^{\f16},\\
|\p_1\td{W}\circ\Psi_W^{y_0}(s)|\leq \f12\ve^\f{1}{20}\eta^{-\f13},\\
|\ck\nabla\td{W}\circ\Psi_W^{y_0}(s)|\leq \f12\ve^\f{1}{20}.
\end{align*}
\end{subequations}
\end{lem}

\proof
The strategy is to perform weighed esitmates under Lagrangian tranjectories regime.

\textbf{Step 1.} General weighed estimate for $l\leq |y|\leq L$

 Let $|\mu|\leq \f13$, $q=\eta^\mu\p^\ga \td{W}$, recall (\ref{DtdW}) for the evolution of $\p^\ga \td{W}$, 
\begin{equation}
(\p_s +D^{(\ga)}_{\td{W}})\p^\ga \td{W}+(\mathcal{V}_W\cdot\nabla) \p^\ga \td{W}=F^{(\ga)}_{\td{W}},
\end{equation}then we have
\begin{equation}
\p_s q+D_q q+\mathcal{V}_W\cdot\nabla q=\eta^\mu F^{(\ga)}_{\td{W}}, 
\end{equation}where 
\begin{align}
D_q&=D^{(\ga)}_{\td{W}}-3\mu+3\mu\eta^{-1}-2\mu\eta^{-1}\left(y_1(\be_\tau W+G_W)+3h^\nu y_\nu|\ck y|^4    \right)\notag\\
&=D^{(\ga)}_{\td{W}}-3\mu+3\mu\eta^{-1}-2\mu D_\eta.
\end{align}
First we estimate $D_\eta$. By (\ref{pW}), (\ref{G}), \textcolor{red}{ (\ref{Gve})} and (\ref{ph}), we have
\begin{align*}
|D_\eta|\leq&\eta^{-1}\left(2|y_1|\eta^{\f16}+|y_1||G_W|+3|h^\nu||y_\nu||\ck y|^4 \right)\\
\leq&2\eta^{-\f13}+\eta^{-1}|y_1|M(e^{-\f s2}+|y_1|e^{-s}+|\ck y|\ve^{\f12})+3\eta^{-\f16}M\ve^{\f12}e^{-\f s2}\\
\leq&5\eta^{-\f13}+e^{-\f s3},
\end{align*}
as long as $\ve$ is taken small enough.

When $|y_0|\geq l$, since $s_0\geq\log\ve$, $l\ll1$,  by (\ref{int1}),  $D_\eta$ can be controlled as follows
\begin{equation}\label{lint}
2\mu\int_{s_0}^{s}|D_\eta\circ\Psi_W^{y_0}(s')|ds'\leq\int_{s_0}^{\infty}10\left(1+l^2e^{\f25(s'-s_0)} \right)^{-\f13}+e^{-\f{s'}{3}}ds'\leq60\log\f1l+\ve^{\f13}\leq70\log\f1l,
\end{equation}
for all $|\mu|\geq\f12$. Hence by Gr\"onwall's inequality, we have
\begin{align} \label{leqL}
|q\circ\Psi^{y_0}_W(s)|\leq l^{-70}&|q(y_0)|\exp\left(\int_{s_0}^{s}(3\mu-D^{(\ga)}_{\td{W}}-3\mu\eta^{-1})\circ\Psi_W^{y_0}(s')ds' \right)\notag\\
&+l^{-70}\int_{s_0}^{s}|\eta^\mu \td{F}^{(\ga)}_W\circ\Psi_W^{y_0}(s')|
\exp\left(\int_{s'}^{s}(3\mu-D^{(\ga)}_{\td{W}}-3\mu\eta^{-1})\circ\Psi_W^{y_0}(s'')ds'' \right)ds'.
\end{align}

\textbf{Step 2.} Estimate of $\td{W}(y,s)$ for $l\leq|y|\leq L$.

Taking $\mu=-\f16$ then the damping exponent is $3\mu-D^{(\ga)}_{\td{W}}-3\mu\eta^{-1}=-\be_\tau\p_1\bar{W}+\f12\eta^{-1}$. By (\ref{pbarW}) and  (\ref{int1}), we have
\begin{equation}
\int_{s_0}^{s}\left(\be_\tau|\p_1\bar{W}|+\f12\eta^{-1}\right)\circ\Psi^{y_0}_W(s')ds'\leq 5\int_{s_0}^{s}\eta^{-\f13}\circ\Psi^{y_0}_W(s')ds'\leq 30\log\f1l
\end{equation}for all $s\geq s_0\geq-\log\ve$. Next, we apply (\ref{ptdFW}) with $|\ga|=0$ on the forcing term and obtain 
\begin{equation}
\int_{s_0}^{s}|\eta^{-\f16}\td{F}_W|\circ\Psi^{y_0}_W(s')ds'\les \ve^{\f{1}{10}}\int_{s_0}^{s}\eta^{-\f16}\circ\Psi^{y_0}_W(s')ds'\les \ve^{\f{1}{10}}\log\f1l
\end{equation}
for all $s\geq s_0\geq-\log\ve.$ Plugging the above two bounds into (\ref{leqL}), we have
\begin{equation}
|\eta^{-\f16}\td{W}\circ\Psi^{y_0}_W(s)|\leq l^{-100}\eta^{-\f16}(y_0)|\td{W}(y_0, s_0)|+M\ve^{\f{1}{10}}l^{-100}\log\f1l,
\end{equation}
where  the implicit constant of force estimate is absorbed in $M$. Then we use initial data (\ref{initdW}) at $s_0=-\log\ve$ and bootstrap assumption (\ref{ptdW3}) for $s\geq-\log\ve$ to get
\begin{equation}
\eta^{-\f16}(y)|\td{W}(y,s)|\leq l^{-100}\max\{M\ve^{\f{1}{18}}l^4,\ve^\f{1}{4} \}+M\ve^{\f{1}{15}}l^{-100}\log\f1l\leq \f{1}{10}\ve^\f{1}{20}.
\end{equation}

\textbf{Step 3.} Estimate of $\p_1\td{W}(y,s)$ for $l\leq|y|\leq L$.

We adapt the same argument as in \textbf{Step 2.}. First, we take $\mu=\f13$ and the damping exponent is $3\mu-D^{(1,0,0)}_{\td{W}}-3\mu\eta^{-1}=-\be_\tau(\p_1W+\p_1\bar{W})-\eta^{-1}$. Then it holds that
\begin{equation}
\int_{s_0}^{s}\be_\tau(|\p_1W+\p_1\bar{W}|+\eta^{-1})\circ\Psi^{y_0}_W(s')ds'\leq80\log\f1l.
\end{equation}
Due to (\ref{pFW2}), the forcing term is bounded by
\begin{equation}
\int_{s_0}^{s}|\eta^{\f13}\td{F}^{(1,0,0)}_W|\circ\Psi^{y_0}_W(s')ds'\les \ve^\f{1}{15}\log\f1l.
\end{equation}
Therefore, \eqref{initdW}, combined with the above two estimate, gives 
\begin{equation*}
\eta^{\f13}(y)|\p_1\td{W}(y,s)|\leq l^{-140}\max\{M\ve^{\f{1}{18}}l^3,\ve^{\f{1}{4}} \}+M\ve^{\f{1}{15}}l^{-140}\log\f1l\leq \f{1}{10}\ve^{\f{1}{20}}.
\end{equation*}

\textbf{Step 4.} Estimate of $\ck\nabla\td{W}(y,s)$ for $l\leq |y|\leq L$.

Without loss of generality we assume $\ga=(0,1,0)$. We take $\mu=0$ and the damping exponent is 
$3\mu-D^{(0,1,0)}_{\td{W}}-3\mu\eta^{-1}=-\be_\tau\p_1\bar{W}$. Then from (\ref{pFW2}), (\ref{int2}), we have
$$|\ck\nabla{\td{W}}(y,s)|\leq l^{-100}\max\{M\ve^{\f{1}{18}}l^3, \ve^{\f14} \}+M\ve^\f{1}{15}\leq \f{1}{10}\ve^{\f{1}{20}}.$$

\textbf{Step 5.} Estimate of $\td W(y,s)$, $\nabla\td W(y,s)$ for $|y|\leq l$.

When $|y|\leq l$ the bounds can be obtained by Taylor expansion of $\td W$ and $\nabla\td W$ at $y=0$ with $\td W(0,s)=\nabla\td W(0,s)=0$.

\qed

\subsection{Closure of $L^\infty$ Type Bootstrap on $W$}

\begin{lem}
If set $\Psi=\Psi_W^{y_0}$, we have
\begin{subequations}
\begin{align}
|\eta^{-\f16}W\circ\Psi(s)|\leq&\f32,\label{636a}\\
|\eta^{\f13}\p_1W\circ\Psi(s)|\leq&\f32,\label{636b}\\
|\ck\nabla W\circ\Psi(s)|\leq&\f32,\label{636c}\\
|\eta^{\f13}\p^\ga W\circ\Psi(s)|\leq& \f32M^{\f{1+|\ck\ga|}{3}},\ |\ga|=2, \ga_1\geq1,\label{636d}\\
|\eta^{\f16}\ck\nabla^2W\circ\Psi(s)|\leq&\f32M.\label{636e}
\end{align}
\end{subequations}
\end{lem}
\proof
Recall (\ref{Devo1}), 
\begin{equation}
\left( \p_s+D^{(\ga)}_W \right) \p^\ga W+(\mathcal{V}_W\cdot\nabla)\p^\ga W=F^{(\ga)}_W,
\end{equation}where
$D^{(\ga)}_W=\f{3\ga_1+\ga_2+\ga_3-1}{2}+\be_\tau\ga_1\p_1W$.

Let $q=\eta^\mu\p^\ga W$, then when \textcolor{red}{$y_0$?}$l\leq|y|\leq L$, we have
\begin{align} \label{Wgeql}
|q\circ\Psi^{y_0}_W(s)|\leq Ml^{-70}&|q(y_0)|\exp\left(\int_{s_0}^{s}(3\mu-D^{(\ga)}_W-3\mu\eta^{-1})\circ\Psi_W^{y_0}(s')ds' \right)\notag\\
&+Ml^{-70}\int_{s_0}^{s}|\eta^\mu F^{(\ga)}_W \circ\Psi_W^{y_0}(s')|
\exp\left(\int_{s'}^{s}(3\mu-D^{(\ga)}_W-3\mu\eta^{-1})\circ\Psi_W^{y_0}(s'')ds'' \right)ds'.
\end{align}
When $|y_0|\geq L$, by (\ref{int2}) with $\de=\f13$, we get
\begin{equation}
2\mu\int_{s_0}^{s}|D_\eta\circ\Psi_W^{y_0}(s')|ds'\leq\int_{s_0}^{\infty}10\left(1+L^2e^{\f25(s'-s_0)}\right)^{-\f13}+e^{-\f{s'}{3}}ds'\les L^{-\f23}+\ve^{\f13}\les \ve^{\f{1}{16}},
\end{equation}
which yields
\begin{align} \label{WgeqL}
|q\circ\Psi^{y_0}_W(s)|\leq e^{\ve^{\f{1}{16}}}&|q(y_0)|\exp\left(\int_{s_0}^{s}(3\mu-D^{(\ga)}_{W}-3\mu\eta^{-1})\circ\Psi_W^{y_0}(s')ds' \right)\notag\\
&+ e^{\ve^{\f{1}{16}}}\int_{s_0}^{s}|\eta^\mu F^{(\ga)}_W\circ\Psi_W^{y_0}(s')|
\exp\left(\int_{s'}^{s}(3\mu-D^{(\ga)}_{W}-3\mu\eta^{-1})\circ\Psi_W^{y_0}(s'')ds'' \right)ds'.
\end{align}

Now we note that \eqref{636a}-\eqref{636c} hold due to (\ref{pbarW}) and $W=\td{W}+\bar{W}$. So we only need to prove \eqref{636a}-\eqref{636c} for $|y|\geq L$ and \eqref{636d}-\eqref{636e} for $|y|\geq l$.

\textbf{Step 1.} Estimate of $W(y,s)$ for $|y|\geq L$.

We take $\mu=-\f16$, then the damping term is $3\mu-D^{(\ga)}_W-3\mu\eta=\f12\eta^{-1}$ and the forcing term is $\eta^{-\f16}(F_W-e^{-\f s2}\be_\tau\dot{\ka})$. By (\ref{int2}) with $\de=1$, it holds that $$\int_{s_0}^{s}\f12\eta^{-1}\circ\Psi^{y_0}_W(s')ds'\leq\int_{s_0}^{\infty}(1+L^2e^{\f25(s'-s_0)})^{-1}ds'\les L^{-1}\les\ve^{\f{1}{10}},$$
as well as $$\int_{s_0}^{s}|(F_W-e^{-\f s2}\be_\tau\dot{\ka})\circ\Psi^{y_0}_W(s')|\les\int_{s_0}^{s}e^{-\f {s'}{2}}\les \ve^\f12.$$

Then (\ref{WgeqL}) implies that 
$$|\eta^{-\f16}W\circ\Psi^{y_0}_W(s)|\leq e^{2\ve^{\f{1}{16}}}(|\eta^{-\f16}W(y_0,s_0)|+\ve^{\f13}).$$

When $s_0>-\log\ve$, $|y_0|=L$,  (\ref{ptdW1}) gives $|\eta^{-\f16}W(y,0)|\leq 1+\ve^{\f{1}{20}}$. If $s_0=-\log\ve$, then by initial assumption (\ref{initdW}), there still holds $|\eta^{-\f16}W(y,0)|\leq 1+\ve^{\f{1}{20}}$. In conclusion, for any $|y|\geq L$, 
$$|\eta^{-\f16}W\circ\Psi^{y_0}_W(s)|\leq e^{2\ve^{\f{1}{16}}}(1+\ve^{\f{1}{11}}+\ve^{\f13})\leq \f32.$$

\textbf{Step 2.} Estimate of $\p_1W(y,s)$ for $|y|\geq L$.

We take $\mu=\f13$, then the damping term is $|3\mu-D^{(\ga)}_W-3\mu\eta^{-1}|\leq |\be_\tau\p_1W| +|\eta^{-1}|\leq 3\eta^{-\f13}.$
For $|y_0|\geq L$, by (\ref{int2}) with $\de=\f13$ we have
$$\int_{s_0}^{s}(3\mu-D^{(\ga)}_W-3\mu\eta^{-1})\circ\Psi^{y_0}_W(s')ds'\leq 3\int_{s_0}^{s}\eta^{-\f13}\circ\Psi^{y_0}_W(s')ds'\les L^{-\f23}\les \ve^{\f{1}{15}}.$$
Meanwhile, (\ref{pFW2}) yields
$$\int_{s_0}^{s}|\eta^{\f13}F^{(1,0,0)}_W\circ\Psi^{y_0}_W(s')|ds'\les\ve^{\f18}\int_{s_0}^{s}\eta^{\f13-\f12+\f{1}{3(2k-5)}}\les \ve^\f18\int_{s_0}^{s}\eta^{\f{1}{15}}\circ\Psi^{y_0}_W(s')\les \ve^{\f18}$$
for $|y_0|\geq L.$
Gather these estimate, by (\ref{WgeqL}), we have 
$$|\eta^{\f13}\p_1W\circ\Psi^{y_0}_W(s)|\leq e^{2\ve^{\f{1}{16}}}(|\eta^{-\f13}\p_1 W(y_0,s_0)|+\ve^{\f17}).$$
When $s_0>-\log\ve$ we let $|y_0|=L$ and then
$|\eta^{\f13}\p_1W|\leq\eta^{\f13}|\p_1\bar{W}|+\eta^{\f13}|\p_1\td{W}|\leq 1+\ve^{\f{1}{20}}$. When $s_0=-\log\ve$, since we have $|y_0|\geq L$ and the initial data assumption, we get $|\eta^{\f13}\p_1W|\leq 1+\ve^{\f{1}{20}}.$ Therefore, we reach the conclusion that $$|\eta^{\f13}\p_1W(y,s)|\leq e^{2\ve^\f{1}{16}}(1+\ve^{\f{1}{20}}+\ve^\f17)\leq \f32$$
for all $|y|\geq L$ and $s\geq -\log\ve.$

\textbf{Step 3.} $\ck\nabla W(y,s)$ for $|y|\geq L$ 

We take $\mu=0$ and the damping term is given by $3\mu-D^{(\ga)}_W-3\mu\eta^{-1}=-\be_\tau\p_1 W$, which we will appeal same estimate as in previous proof.
For the forcing term we have 
$$\int_{s_0}^{s}|F^{(0,1,0)}_W\circ\Psi^{y_0}_W(s')|ds'\leq M\ve^{\f13}\eta^{-\f13}\leq \ve^\f14.$$
For $s_0>-\log\ve$ we still use $W=\td{W}+\bar{W}$ while for $s_0=-\log\ve$ we use initial data assumption. So at last  
$$|\ck\nabla W(y,s)|\leq e^{2\ve^{\f{1}{16}}}(1+\ve^\f14)\leq \f32$$
holds for all $|y|\geq L$ and $s\geq -\log\ve.$

\textbf{Step 4.} $\p^\ga W(y,s)$ with $|\ga|=2$ for $|y|\geq l$

We take $\mu= \f13$ for $|\ga|=2,\ \ga_1\geq1$ and $\mu=\f16$ for $|\ga|=2$, $\ga_1=0$, respectively. 

Case 1. $\ga_1=0$, $|\ck\ga|=2$.
The damping term becomes $3\mu-D^{(\ga)}_W-3\mu\eta^{-1}=-\be_\tau\p_1 W$, so \begin{equation}\int_{s_0}^{s}\be_\tau |\p_1 W|\circ\Psi^{y_0}_W(s')ds'\leq 40\log\f1l,\label{661}\end{equation} and forcing term can be bounded as 
\begin{equation}\int_{s_0}^{s}|\eta^{\f16}F^{(\ga)}_W|\circ\Psi^{y_0}_W(s')ds'\leq M^\f23\eta^{-\f16+\f{3}{2k-7}}\leq 6M^\f23\log\f1l.\label{662}\end{equation} 
Then (\ref{Wgeql}), together with \eqref{661} and \eqref{662}, gives
\begin{align*}
\eta^{\f16}|(y)\ck\nabla^2 W(y,s)|\leq& l^{-100}\eta^{\f16}(y_0)|\ck\nabla^2 W(y_0,s_0)|+M^\f23l^{-100}\log\f1l\\
\leq &l^{-100}\max\{ 1,M\ve^{\f{1}{18}}   \}+M^\f23l^{-102}\\
\leq &\f32M,
\end{align*}
where the second inequality we use initial data (\ref{iniW}) for $s_0=-\log\ve$ and bootstrap bound (\ref{pW}) for $s_0>-\log\ve.$

Case 2. $\ga_1\geq1$, $|\ga|=2$. The damping term $3\mu-D^{(\ga)}_W-3\mu\eta^{-1}=-\f{2\ga_1-1}{2}-(2\ga_1-1)\be_\tau\p_1 W-\eta^{-1}$, 
whose exponential integral is bounded by
$$\exp \left(  \int_{s'}^{s}(3\mu-D^{\ga)}_W\circ\Psi^{y_0}_W(s'')ds'' )    \right) \leq l^{-120}e^{\f{s'-s}{2}}.$$
On the other hand the forcing term
$$\int_{s_0}^{s}\eta^{\f13}F^{(\ga)}_We^{\f{s'-s}{2}}ds'\leq2M^{\f{|\ck\ga|}{3}}.$$
So at last 
\begin{align*}
\eta^{\f13}(y)|\p^\ga W(y,s)|\leq& l^{-180}\eta^\f13(y_0)|\p^\ga W(y_0,s_0)|+2M^{\f{|\ck\ga|}{3}}l^{-180}\\
\leq&l^{-180}\max\{1, 2M\ve^{\f{1}{18}}l^2\}+2M^{\f{|\ck\ga|}{3}}l^{-180}\\
\leq&\f32M^{\f{|\ck\ga|+1}{3}}.
\end{align*}

\textbf{Step 5.} $W(y,s)$, $\nabla W(y,s)$ for $|y|\leq L$ and $\nabla^2W(y,s)$ for $|y|\leq l$

The first two bounds can be obtained by the sum of (\ref{ptdWW}) and (\ref{pbarW}), while the last results from Taylor expansion of $\nabla^2W(y,s)$ at $y=0$, in view of $\nabla^2\td W(0,s)=\nabla^2 W(0,s)=0$ and (\ref{ptdW3}).

\qed

\section{Energy Estimates}
In this section we will perform high order energy estimates to the system about $U_i$ and $S$. We use the homogeneous Sobolev space instead of inhomogeneous one because $\|U\|_{L^2}$, $\|S\|_{L^2}$ can not be taken enough small at initial data.

Recall from (\ref{DevoU}) and (\ref{DevoS}), $U_i$ and $S$ satisfy the equations
$$\p_s U_i+(\mathcal{V}_U\cdot\nabla)U_i=F_{U_i},\ \ \ 
\p_s S+(\mathcal{V}_U\cdot\nabla)S=F_S,$$
where the velocity fields are $$\mathcal{V}_U=(\f32y_1+g_U,\f12y_\nu+h^\nu),$$
and the external forces are
$$F_{U_i}=-2\be_3\be_\tau e^\f s2\de^{i1}  S\p_1S-2\be_1\be_\tau e^\f s2\de^{i1}\p_1\Phi-2\be_3\be_\tau e^{-\f s2}\de^{i\nu}S\p_\nu S-2\be_1\be_\tau e^{-\f s2}\de^{i\nu}\p_\nu \Phi,$$
$$F_S=-2\be_3\be_\tau S(e^\f s2\p_1U_1+e^{-\f s2}\p_\nu U_\nu).$$
Now we apply $\p^\ga$ to both sides of above equations with $|\ga|=k$ to get 
\begin{lem}
For $|\ga|=k$, it holds that
\begin{align}
\p_s(\p^\ga U_i)+(\mathcal{V}_U\cdot\nabla)(\p^\ga U_i)+\mathcal{D}_\ga\p^\ga U_i+2\be_\tau\be_3S[e^{\f s2}\de^{i1}\p_1(\p^\ga S)+e^{-\f s2}\de^{i\nu}\p_\nu(\p^\ga S)]\notag\\
+\be_\tau\be_3(1+\ga_1)\p_1W\p^\ga S=-2\be_\tau\be_1[e^{\f s2}\de^{i1}\p_1(\p^\ga \Phi)+e^{-\f s2}\de^{i\nu}\p_\nu(\p^\ga \Phi)]+\mathcal{F}^{(\ga)}_{U_i},\label{DkU}
\end{align}and 
\begin{align}
\p_s(\p^\ga S)+(\mathcal{V}_U\cdot\nabla)(\p^\ga S)+\mathcal{D}_\ga\p^\ga S+2\be_\tau\be_3S[e^{\f s2}\p_1(\p^\ga U_1)+e^{-\f s2}\p_\nu(\p^\ga U_\nu)]\notag\\
+\be_\tau(\be_1+\be_3\ga_1)\p^\ga U_1\p_1W =\mathcal{F}^{(\ga)}_S \label{DkS}.
\end{align}
Here, the damping term is 
\begin{equation*}
\mathcal{D}_\ga=\ga_1(1+\p_1g_U)+\f12|\ga|.
\end{equation*}
And the external forces terms are 
$$\mathcal F^{(\ga)}_{U_i}=P_1+P_2+P_3+P_4,$$
$$\mathcal F^{(\ga)}_S=Q_1+Q_2+Q_3+Q_4.$$
with 
\begin{align*}
P_1=-\sum_{\substack{|\be|=|\ga|-1 \\ \be_1=\ga_1}}\binom\ga\be\p^{\ga-\be}g_U\p^\be\p_1U_i-\sum_{\substack{|\be|=|\ga|-1}}\binom\ga\be\p^{\ga-\be}h^\nu\p^\be\p_\nu U_i,
\end{align*}

\begin{align*}
P_2=-\sum_{1\leq|\be|\leq|\ga|-2}\binom\ga\be\big(\p^{\ga-\be}g_U\p^\be\p_1U_i + \p^{\ga-\be}h^\nu\p^\ga\p_\nu U_i\big),
\end{align*}

\begin{align*}
P_3=&-2\be_\tau\be_3e^{-\f s2}\de^{i\nu}\p_\nu S\p^\ga S-2\be_\tau\be_3\sum_{|\be|=|\ga|-1}\binom\ga\be e^{-\f s2}\de^{i\nu}\p^{\ga-\be}S\p^\be\p_\nu S\\
&+\be_\tau\be_3e^{\f s2}(1+\ga_1)\de^{i1}\p_1Z\p^\ga S-2\be_\tau\be_3\sum_{|\be|=|\ga|-1}\binom\ga\be e^{\f s2}\p^{\ga-\be}S\p^\ga\p_1S,
\end{align*}

\begin{align*}
P_4=-2\be_\tau\be_3\sum_{1\leq|\be|\leq|\ga|-2}\binom\ga\be\big(e^{\f s2}\p^{\ga-\be}S\p^\be\p_1S+e^{-\f s2}\de^{i\nu}\p^{\ga-\be}S\p^\be\p_\nu S \big),
\end{align*}
and 

\begin{align*}
Q_1-\sum_{\substack{|\be|=|\ga|-1\\ \be_1=\ga_1}}\binom\ga\be\p^{\ga-\be}g_U\p^\be\p_1 S-\sum_{|\be|=|\ga|-1}\binom\ga\be\p^{\ga-\be}h^\nu\p^\be\p_\nu S,
\end{align*}

\begin{align*}
Q_2=&-2\be_\tau\be_1e^{-\f s2}\p_\nu S\p^\ga U_\nu-2\be_\tau\be_3\sum_{|\be|=|\ga|-1}\binom\ga\be e^{-\f s2}\p^{\ga-\be}S\p^\be\p_\nu U_\nu\\
&+\be_\tau(\be_1+\be_3\ga_1)e^{\f s2}\p_1Z\p^\ga U_1-2\be_\tau\be_3\sum_{\substack{|\be|=|\ga|-1\\ \be_1=\ga_1}}\binom\ga\be e^{\f s2}\p^{\ga-\be}S\p^\be\p_1U_1,
\end{align*}

\begin{align*}
Q_3=&-\sum_{1\leq|\be|\leq|\ga|-2}\binom\ga\be\big(\p^{\ga-\be}g_U\p^\be\p_1S+\p^{\ga-\be}h^\nu\p^\be\p_\nu S \big)   \\
&-2\be_\tau\be_3\sum_{1\leq|\be|\leq|\ga|-2}\binom\ga\be \big(e^{\f s2}\p^{\ga-\be}S\p^\be\p_1U_1+e^{-\f s2}\p^{\ga-\be}S\p^\ga\p_\nu U_\nu \big).
\end{align*}
\end{lem}

\proof
(\ref{DkU}) follows from direct computation, while the (\ref{DkS}) is similar. 

The main focus is on the third term of $P_3$. Since $\p_1S=\f12(e^{-\f s2}\p_1W-\p_1Z)$, we have
\begin{align*}
&\p^\ga(-2\be_3\be_\tau e^{\f s2}\de^{i1}S\p_1S)\\
=&-2\be_3\be_\tau e^{\f s2}\de^{i1}S\p_1\p^\ga S-2\be_3\be_\tau e^{\f s2}\de^{i1}\p^\ga S\p_1S-2\be_3\be_\tau e^{\f s2}\de^{i1}\ga_1\p_1S\p^\ga S\\
&-2\be_3\be_\tau e^{\f s2}\de^{i1}\sum_{|\be|=|\ga|-1}\binom\ga\be\p^{\ga-\be}S \p^\be\p_1S
-2\be_3\be_\tau\sum_{1\leq|\be|\leq|\ga|-2}e^{\f s2}\de^{i1}\binom\ga\be\p^{\ga-\be} S\p^\be\p_1S,
\end{align*}
where the second and the third term on the right hand side can be further rewriten as 
\begin{align*}
-2\be_3\be_\tau e^{\f s2}\de^{i1}\p^\ga S\p_1S-2\be_3\be_\tau e^{\f s2}\ga_1\p_1S\p^\ga S=&-2\be_\tau\be_3(1+\ga_1)e^{\f s2}\de^{i1}\p^\ga S\p_1S\\
=&-\be_\tau\be_3(1+\ga_1)\de^{i1}\p^\ga S(\p_1W+e^{\f s2}\p_1Z).
\end{align*}
Thus (\ref{DkU}) follows.

\qed

Now we define the energy $E_k(s)$ by refined $\dot{H}^k$ norm for $k\geq18$
\begin{equation}
E^2_k(s)=\sum_{|\ga|=k}\la^{|\ck\ga|} \left(\|\p^\ga U(\cdot,s)\|^2 +\|\p^\ga S(\cdot,s)\|^2 \right),
\end{equation}
where $\la=\la(k)=\f{\de^2}{12k^2}$ and $\de\in(0,\f{1}{32}]$ is to be determined. We can easily see that $E_k(s)$ is equivalent to $\dot H^k$ norm despite the factor $\la^{|\ck\ga|}$. The introduction of the factor is just a technical issue.

\begin{lem}
For $k\geq18$ and $|\ga|=k$, the estimate for external force $\mathcal{F}^{(\ga)}_{U_i}$ and $\mathcal{F}^{(\ga)}_{S}$ are
\begin{subequations}
\begin{align*}
2\sum_{|\ga|=k}\la^{|\ck\ga|}\int_{\mathbb{R}^3}|\mathcal{F}^{(\ga)}_{U_i}\p^\ga U_i|\les (2+8\de)E^2_k+e^{-s}M^{4k-1},\\
2\sum_{|\ga|=k}\la^{|\ck\ga|}\int_{\mathbb{R}^3}|\mathcal{F}^{(\ga)}_{S}\p^\ga S|\les (2+8\de)E^2_k+e^{-s}M^{4k-1}.
\end{align*}
\end{subequations}
\end{lem}
\proof 
The proof is fairly similar to that of the Lemma 12.2 in \cite{Buckmaster2019} and we omit it.

\qed

\begin{lem}
For the electric potential terms, we have 
\begin{subequations}
\begin{align}
\int e^{\f s2}(\p^\ga U_i)\p_1(\p^\ga\Phi)\leq e^{-s}\|U\|_{\dot{H}^k}(\|S\|_{\dot{H}^k}+\|R_+\|_{\dot H^k}),\label{p1Phi}\\
\int e^{-\f s2}(\p^\ga U_i)\p_\nu(\p^\ga\Phi)\leq e^{-s}\|U\|_{\dot{H}^k}(\|S\|_{\dot{H}^k}+\|R_+\|_{\dot H^k}).\label{pnuPhi}
\end{align}
\end{subequations}
\end{lem}
\proof
By H\"older inequality, we have $$
\int e^{\f s2}(\p^\ga U_i)\p_1(\p^\ga\Phi)\leq e^{\f s2 }\|\p^\ga U\|_{L^2}\|\p^\ga\p_1\Phi\|_{L^2}.$$
Take $\p_1$ to (\ref{poten}), we get$$\p_1\Phi=\int_{\mathbb{R}^3}\f{-e^{-\f52s} e^{-3s}z_1(R_+-R_e)(y-z)  }{(e^{-3s}|z_1|^2+e^{-s}|\ck z|^2  )^\f32},$$
 then  (\ref{P1}) implies that
\begin{align*}
\|\p^\ga\p_1\Phi\|_{L^2}=e^{-\f32s}\|P_1\star\p^\ga (R_+-R_e)\|_{L^2}\les e^{-\f32s}(\|R_e\|_{\dot{H}^k}+\|R_+\|_{\dot H^k}).
\end{align*}
At last, it suffices to show that $$\|R_e\|_{\dot{H}^k}\les \|\p^\ga S\|_{L^2}\les \|S\|_{\dot{H}^k}.$$
 Indeed, for $|\ga|=k$, write $\p^\ga (S^\f1\al)$ by using
 Fa\'a di Bruno formula (see \cite{CS96}), for $|\ga|=k$,
\begin{align*}
|\p^\ga (S^\f1\al)|\les \sum_{|\ga|=1}^k S^{\f1\al-|\nu|}\sum_{P(\be,\nu)}
\prod_{l=1}^k
\f{|\p^{\be_l}S|^{\nu_l}}{\nu_l!(\be_l!)^{\nu_l}},
\end{align*}
Then we taking $\f12=\sum_l\f{1}{q_l}$ and $\sum_{l}\be_l\nu_l=k$, $\sum_{l}\nu_l=|\nu|$,  by H\"older inequality and Sobolev interpolation,
\begin{align*}
\sum_{|\ga|=k}\|\p^\ga (S)^\f1\al\|_{L^2} 
\les&\prod_{l=1}^k\|S^{\f1\al-|\nu|}\|_{L^\infty}\|(\p^{\be_l} S)^{\nu_l}\|_{L^{q_l}}\\
\les&\prod_{l=1}^k\|S\|^{\f1\al-|\nu|}_{L^\infty}\|\p^{\be_l}S\|^{\nu_l}_{L^{q_l\nu_l}}\\
\les&\prod_{l=1}^k\|S\|^{\f1\al-|\nu|}_{L^\infty}\|\p^{\be_l}S\|^{\nu_l(1-\f{\be_l}{k})}_{L^{\infty}}\|S\|^{\nu_l\f{\be_l}{k}}_{\dot H^k}\\
\les&\|S\|^{\f1\al-1}_{L^\infty}\|S\|_{\dot H^k}.
\end{align*}

\qed

\begin{prop}
For $k\geq 18$, we have
\begin{equation}\label{Hk2}
E^2_k(s)\leq e^{-2s}\ve^{-1}+2e^{-s}M^{4k-1}(1-\ve^{-1}e^{-s}).
\end{equation}
\end{prop}
\proof

Taking $L^2$ inner product on (\ref{DkU}) and (\ref{DkS}) with $\la^{|\ck\ga|}\p^\ga U_i$ and $\la^{|\ck\ga|}\p^\ga S$, respectively. Then we use integration by parts:
\begin{align*}
\p_s\int\la^{|\ck\ga|}(|\p^\ga U|^2+|\p^\ga S|^2)+\la^{|\ck\ga|}\int(2\mathcal{D}_\ga-\nabla\cdot\mathcal{V}_U)(|\p^\ga U|^2+|\p^\ga S|^2)\\
+2\be_\tau\la^{|\ck\ga|}\int_{\mathbb{R}^3}(\be_1+\be_3+2\be_3\ga_1)\p_1W\p^\ga S\p^\ga U_1    \\
=4\be_\tau\be_3\la^{|\ck\ga|}\int\left[e^{\f s2}(\p_1S)(\p^\ga U_1)+e^{-\f s2}(\p_\nu S)(\p^\ga U_\nu)\right](\p^\ga S)\\
-4\be_\tau\be_3\la^{|\ck\ga|}\int\left[e^{\f s2}(\p^\ga U_1)\p_1(\p^\ga \Phi)+e^{-\f s2}(\p^\ga U_\nu)\p_\nu(\p^\ga \Phi)\right]\\
+2\la^{|\ck\ga|}\int \mathcal{F}^{(\ga)}_{U_i}\p^\ga U_i+\mathcal{F}^{(\ga)}_S \p^\ga S.
\end{align*}
Then summing over all $|\ga|=k$ we get 
\begin{equation}\label{psE}
\p_s E^2_k(s)+A_1+A_2=B_1+B_2+B_3.
\end{equation}
We first estimate the damping term $A_1$ as
\begin{align*}
&2\mathcal{D}_\ga-\nabla\cdot\mathcal{V}_U\\
=&2\ga_1(1+\p_1g_U)+|\ga|-\f52-(\p_1 g_U+\p_2h^2+\p_3h^3)\\
=&(2\ga_1-1)(\be_\tau\be_1\p_1W+\p_1G_U)+2\ga_1+|\ga|-\f52-\p_2h^2-\p_3h^3\\
\geq&|\ga|-\f52+2\ga_1-\be_\tau\be_1(1+\ve^\f{1}{20})(2\ga_1-1)-\ve^\f15,
\end{align*}where we have used (\ref{pWlbd}), (\ref{pG}) (\ref{ph}) and $\be_1+\be_3=1$. 
 $A_2$ also has lower bound as following,
\begin{align*}
&2\be_\tau\la^{|\ck\ga|}(\be_1+\be_3+2\be_3\ga_1)\p_1W\p^\ga S\p^\ga U_1\\
\geq&-\be_\tau(\be_1+\be_3+2\be_3\ga_1)|\p_1W|(|\p^\ga S|^2+|\p^\ga U|^2)\\
\geq&-\be_\tau(1+2\be_3\ga_1)(|\p^\ga S|^2+|\p^\ga U|^2).
\end{align*}
So $A_1+A_2$ is bounded from below as
\begin{align*}
\sum_{|\ga|=k}\left[|\ga|-\f52+2\ga_1-\be_\tau\be_1(1+\ve^\f{1}{20})(2\ga_1-1)-\ve^\f15-\be_\tau(1+2\be_3\ga_1)\right]\int_{\mathbb{R}^3}\la^{|\ck\ga|}(|\p^\ga S|^2+|\p^\ga U|^2)\\
\geq\sum_{|\ga|=k}\left[|\ga|-\f52+2\ga_1(1-\be_\tau)+\be_\tau\be_1-\be_\tau-\ve^\f{1}{30}  \right]\int_{\mathbb{R}^3}\la^{|\ck\ga|}(|\p^\ga S|^2+|\p^\ga U|^2)\\
\geq (k-5)\sum_{|\ga|=k}\int_{\mathbb{R}^3}\la^{|\ck\ga|}(|\p^\ga S|^2+|\p^\ga U|^2).
\end{align*}

For $B_1$, (\ref{pW}) and (\ref{pZ}) give
\begin{align*}
B_1\leq &4\be_\tau\be_3\sum_{|\ga|=k}\la^{|\ck\ga|}(\|\p_1W\|_{L^\infty}+e^\f s2\|\p_1Z\|_{L^\infty}+e^{-s}\|\ck\nabla W\|_{L^\infty}+e^{-\f s2}\|\ck\nabla Z\|_{L^\infty})\| \p^\ga U \|_{L^2}\|\p^\ga S\|_{L^2}   \\
\leq&2(1+M\ve)(1+\ve^\f12)E^2_k(s)\\
\leq&(2+\ve^\f15)E^2_k(s).
\end{align*}
Furthermore, by Lemma 8.1, Lemma 8.2 and (\ref{R+Hk}), we get
$$B_2+B_3\leq 2(2+8\de)E^2_k(s)+2e^{-s}M^{4k-1}+e^{-s} E^2_k(s),$$
where the factor $\ve$ can be absorbed in $M^{4k-1}$.
Therefore, from $k\geq 18$ and $\de\leq \f{1}{32}$ and (\ref{psE}), we finally get 
$$\f{d}{ds}E^2_k(s)+2E^2_k(s)\leq 2e^{-s}M^{4k-1},$$ so 
$$E^2_k(s)\leq e^{-2(s_0-s)}E^2_k(s_0)+2e^{-s}M^{4k-1}(1-e^{-(s-s_0)}).$$

\qed

\proof[Proof of Proposition 4.4] 
Fisrt by (\ref{iniHk}), we have $E^2_k(-\log\ve)\leq \ve$, so
\begin{align*}
\| U\|^2_{\dot H^k}+\|S \|^2_{\dot H^k}\leq& \la^{-k}\left(\ve^{-1}e^{-2s}+2e^{-s}M^{4k-1}(1-\ve^{-1}e^{-s}) \right)\\
\leq &\la^{-k}\ve^{-1}e^{-2s}+2e^{-s}\la^{-k}M^{4k-1}(1-\ve^{-1}e^{-s}).
\end{align*}
Then,
\begin{align*}
\|Z\|^2_{\dot H^k}\leq &2\left(\|U_1\|^2_{\dot H^k}+\|S\|^2_{\dot H^k} \right)\\
\leq&2\la^{-k}e^{-s}+e^{-s}(1-e^{-s}\ve^{-1})M^{4k},
\end{align*}
where we use $\ve^{-1}\leq e^s$ and $2\la^{-k}<M$. And for $\|W\|_{\dot H^k}$ is similar
\begin{align*}
\|W\|^2_{\dot H^k}\leq& 2e^s\left(\|U_1\|^2_{\dot H^k}+\|W\|^2_{\dot H^k} \right)\\
\leq& 2\la^{-k}e^{-s}\ve^{-1}+(1-e^{-s}\ve^{-1})M^{4k}.
\end{align*}
\qed

\section{The Main Theorem in Physical Variable}
For the completion of the article, it necessary to state the main theorem in physical variable. In fact, the initial data set in $(\xcal, \tcal)$ in the following is equivalent to that of in $(y,s)$ assumed in Section 3.

We set the initial time to be $\tcal_0=-\f{2}{1+\al}\ve$, and set 
\begin{equation}\label{inimod8}
\ka_0:=\ka(\tcal_0),\ \ \tau_0:=\tau(\tcal_0),\ \ \xi_0:=\xi(\tcal_0)
\end{equation}
and
$$((u_1)_0(\xcal),(u_2)_0(\xcal),(u_3)_0(\xcal))=u_0(\xcal):=u(\xcal,\tcal_0),\ \ n_{e,0}(\xcal):=n_e(\xcal,\tcal_0),\ \ \si_0:=\f{n_{e,0}^\al}{\al}.$$
We assume the initial data of $n_{e,0}$ satisfies neutrality
\begin{equation}\label{neu}
\int_{\XX_+}(n_+-n_{e,0})=0.
\end{equation}
We also assume that $u_0$ and $n_{e,0}$ satisfy even condition:
\begin{subequations}\label{even1}
\begin{align}
u_0(\xcal_1,\xcal_2,\xcal_3)=u_0(\xcal_1,-\xcal_2,\xcal_3)=u_0(\xcal_1,\xcal_2,-\xcal_3)\\
n_{e,0}(\xcal_1,\xcal_2,\xcal_3)=n_{e,0}(\xcal_1,-\xcal_2,\xcal_3)=n_{e,0}(\xcal_1,\xcal_2,-\xcal_3)
\end{align}
\end{subequations}

We introduce the Riemann type variables
$$\td w_0(\xcal):=(u_1)_0(\xcal)+\si_0(\xcal),\ \ \td z_0(\xcal):=(u_1)_0(\xcal)-\si_0(\xcal),$$ 
and assume the initial data $(\td w_0-\ka_0,\td z_0, (\td u_\nu)_0, n_{e,0})$ is supported in the set (\ref{sptx}),
\begin{equation}\label{sptx8}
\XX=\{|\xcal_1|\leq \ve^\f12, |\ck x|\leq \ve^\f16\}.
\end{equation}

We choose $\td w_0(\xcal)$ such that 
\begin{subequations}\label{slope}
\begin{align}
&\text{the minimum negative slope in the $e_1$ direction }\\
&\p_{\xcal_1}\td w_0 \text{ attains its global minimum at } \xcal=0
\end{align}
\end{subequations}
\begin{equation}\label{ckna}
\td w_0(0)=\ka_0,\ \ \p_{\xcal_1}\td w_0(0)=-\f{1}{\ve},\ \ \ck\nabla_\xcal\td w_0(0)=0,\ \  \nabla^2_\xcal\td w_{0}(0)=0,\ \ 
\nabla^2_\xcal\p_{\xcal_1}\td w_0(0)>0.
\end{equation}

Set 
$$\bar w_\ve(\xcal-\xi(t))=\bar w_\ve(x):=\ve^\f12\bar W(\ve^{-\f32}x_1,\ve^{-\f12}\ck x),$$ and
$$\wt w_0:=\td w_0(\xcal)-\bar w_\ve(\xcal-\xi(t))=w_0(x)-\bar w_\ve(x)=\ve^\f12\td W(y,-\log\ve)+\ka_0.$$

For $\wt w_0$ and $|(\ve^{-\f32}\xcal_1,\ve^{-\f12}\ck\xcal)|\leq 2\ve^{-\f{1}{10}}$, according to (\ref{initdW})
\begin{subequations}
\begin{align}
|\wt w_0(\xcal)-\ka_0|\leq& \ve^{\f{1}{4}}(\ve^3+\xcal_1^2+|\ck\xcal|^6)^{\f16},\\
|\p_{\xcal_1}\wt w_0(\xcal)|\leq& \ve^\f14(\ve^3+\xcal_1^2+|\ck\xcal|^6)^{-\f13},\\
|\ck\nabla_\xcal\wt w_0(\xcal)|\leq& \ve^\f14.
\end{align}
\end{subequations}
Furthermore, for $|(\ve^{-\f32}\xcal_1,\ve^{-\f12}\ck\xcal)|\leq1$, we assume that
\begin{equation}
|\p^\ga_\xcal\wt w_0(\xcal)|\leq\ve^{\f{1}{20}-\f12(3\ga_1+\ga_2+\ga_3)},\ \ \ \ |\ga|=4
\end{equation}
while $\xcal=0$, we assume 
\begin{equation}
|\p^\ga_\xcal\wt w_0(0)|\leq\ve^{\f{1}{20}-\f12(3\ga_1+\ga_2+\ga_3)},\ \ \ \ |\ga|=3
\end{equation}

For $\td w_0$ and $\xcal\in\XX$ but $|(\ve^{-\f32}\xcal_1,\ve^{-\f12}\ck\xcal)|\geq\f12\ve^{-\f{1}{10}}$, we assume that 
\begin{subequations}
\begin{align}
|\td w_0(\xcal)-\ka_0|\leq& (1+\ve^{\f{1}{4}})(\ve^3+\xcal_1^2+|\ck\xcal|^6)^{\f16},\\
|\p_{\xcal_1}\td w_0(\xcal)|\leq& (1+\ve^\f14)(\ve^3+\xcal_1^2+|\ck\xcal|^6)^{-\f13},\\
|\ck\nabla_\xcal \td w_0(\xcal)|\leq&1.
\end{align}
\end{subequations}
For $\xcal\in\XX$, we assume the $\nabla^2\td w_0$ satisfy
\begin{subequations}
\begin{align}
|\p^2_{\xcal_1}\td w_0(\xcal)|\leq& \ve^{-\f32}(\ve^3+\xcal_1^2+|\ck\xcal|^6)^{-\f13},\\
|\p_{\xcal_1}\ck\nabla_\xcal\td w_0(\xcal)|\leq&\ve^{-\f12}(\ve^3+\xcal_1^2+|\ck\xcal|^6)^{-\f13},\\
|\ck\nabla^2_\xcal\td w_0(\xcal)|\leq&(\ve^3+\xcal_1^2+|\ck\xcal|^6)^{-\f16}.
\end{align}
\end{subequations}

For $\td z_0$ and $ (\td u_\nu)_0$, we assume that 
\begin{align}
&|\td z_0(\xcal)|\leq\ve,\ \ \ |\p_{\xcal_1}\td z_0(\xcal)|\leq1,\ \ \ |\ck\nabla\td z_0(\xcal)|\leq\ve^\f12,\notag\\
&|\p^2_{\xcal_1}\td z_0(\xcal)|\leq \ve^{-\f32},\ \ \ |\p_{\xcal_1}\ck\nabla_\xcal\td z_0(\xcal)|\leq \ve^{-\f12},\ \ \ |\ck\nabla^2_\xcal\td z_0(\xcal)|\leq 1,
\end{align}and 
\begin{equation}
|(\td u_\nu)_0(\xcal)|\leq\ve,\ \ \ |\p_{\xcal_1}(\td u_\nu)_0(\xcal)|\leq1,\ \ \ 
|\ck\nabla_\xcal(\td u_\nu)_0(\xcal)|\leq\ve^\f12,\ \ \ |\ck\nabla^2_\xcal(\td u_\nu)_0(\xcal)|\leq 1.
\end{equation}

For the initial specific vorticity, we assume that
\begin{equation}\label{Hk8}
\left\| \f{\nabla_\xcal\times u_0(\xcal)}{n_{e,0}} \right\|_{L^\infty}\leq1.
\end{equation}

For the Sobolev norm of initial condition we assume for a fixed $k\geq18$
\begin{equation} 
\sum_{|\ga|=k}\ve^2\|\p^\ga_\xcal\td w_0\|^2_{L^2}+\|\p^\ga_\xcal\td z_0\|^2_{L^2}+\|\p^\ga_\xcal(\td u_\nu)_0\|^2_{L^2}\leq \ve^{\f72-(3\ga_1+|\ck\ga|)}.
\end{equation}

At last recall the transform
\begin{align}\label{tran8}
&(x,t)=(\xcal-\xi(\f{1+\al}{2}\tcal),\f{1+\al}{2}\tcal),\ \ \ \td u(x,t)=u(\xcal,\tcal),\ \ \ \td\ze(x,t)=\ze(\xcal,\tcal)=\f{\nabla\times u}{n_e},\notag\\ 
&\td \si(x,t)=\si(\xcal,\tcal)=\f{n_e^\al(\xcal,\tcal)}{\al}=\f{\td n_e^\al(x,t)}{\al}.
\end{align}

\begin{thm}
Let $\ga>1$, $\al=\f{\ga-1}{2}$, $T_*>0$. $u_0$, $n_{e,0}$, $\si_0$, $\td w_0$, $\td z_0$ are defined above, where $n_{e,0}$ satisfying the neutrality \emph{(\ref{neu})} and even condition \emph{(\ref{even1})}. The modulation variables $\ka$, $\tau$, $\xi$ have initial conditions given by \emph{(\ref{inimod8})}. We also assume the initial data $(\td w_0-\ka_0,\td z_0, (\td u_\nu)_0, n_{e,0})$ is supported in the set \emph{(\ref{sptx8})} and satisfy the condition \emph{(\ref{slope})-(\ref{Hk8})}.

 Then there exist $\ka_0>1$, $\ve\ll1$, $T_*=O(\ve)$ and a unique solution 
$$(u,n_e)\in C\left([-\f{2}{1+\al}\ve,-\f{2}{1+\al}T_* );H^k\right)\cap C^1\left([-\f{2}{1+\al}\ve,-\f{2}{1+\al}T_* );H^{k-1}\right)$$ to \emph{(\ref{EP})} which blows up in an asymptotically self-similar type at time $T_*$, at the single point $\xi_*\in\mathbb R^3$.

Also, if use the variable $(x,t)$, $\td u(x,t)$, $\td \si(x,t)$ defined in \emph{(\ref{tran8})}, the following result holds:
\begin{itemize}
\item $T_*=O(\ve^2)$, $\xi_*=O(\ve)$, $|\ka_*-\ka_0|=O(\ve)$, where $T_*$ is defined by $\int_{-\ve}^{T_*}(1-\dot\tau(t))dt=\ve$, $\xi_*=\lim_{t\to T_*}\xi(t)$, $\ka_*=\lim_{t\to T_*}\ka(t).$
\item We have $\sup_{t\in[-\ve,T_*)}(\|\td u_1-\f12\ka_0\|_{L^\infty}+\|\td u_\nu\|_{L^\infty}+\|\td\si -\f12\ka_0\|_{L^\infty})\les1.$
\item There holds $\lim_{t\to T_*}\p_{x_1}\td w(\xi(t),t)=-\infty$ and 
$\f{1}{2(T_*-t)}\leq\|\p_{x_1}\td w(\xi(t),t)\|_{L^\infty}\leq \f{2}{T_*-t}$ as $t\to T_*.$
\item The only blowup first order derivatives are $\p_{x_1}\td u$ and $\p_{x_1}\td n_e$, the other first derivatives remain bounded uniformly in $t$:
\begin{subequations}
\begin{align}
&\lim_{t\to T_*}\p_{x_1}u_1(\xi(t),t)=\lim_{t\to T_*}\p_{x_1}\td n_e(\xi(t),t)=-\infty,   \\
&\sup_{t\in[-\ve,T_*)}\|\ck\nabla\td n_e(\cdot,t)\|_{L^\infty}+\|\ck\nabla\td u_1(\cdot,t)\|_{L^\infty}+\|\nabla u_\nu(\cdot,t)\|_{L^\infty}\les1.
\end{align}
\item The electron density is uniformly bounded from below in the support set of $n_e$, especially in the set $\XX$:
$$\sup_{t\in[-\ve,T_*)}\|\td n_{e}^\al(\cdot,t)-\f{\al}{2}\ka_0\|_{L^\infty_\XX}\leq \ve^{\f{1}{10}}. $$
\end{subequations}
\proof

The proof is easy to obtained from bootstrap argument performed in above sections, which we left it to readers.

\qed

\begin{rek}
The assumption \eqref{even1} is made to avoid the over repeanting parts with \cite{Buckmaster2019}. In fact, similar conclusions to \cite{Buckmaster2019} can be obtained if such assumption is removed.

\end{rek}

\end{itemize}
\end{thm}


\begin{thebibliography}{10}

\bibitem{Buckmaster2019}
Tristan Buckmaster, Steve Shkoller, and Vlad Vicol.
\newblock Formation of point shocks for 3d compressible Euler.
\newblock {\em arXiv:1912.04429}, 2019.


\bibitem{Buckmaster2019a}
Tristan Buckmaster, Sammer Iyer.
\newblock Formation of unstable shocks for 2D isentropic compressible
Euler.
\newblock {\em arXiv:2007.15519}, 2020.

\bibitem{Buckmaster2020a}
Tristan Buckmaster, Steve Shkoller, and Vlad Vicol.
\newblock Shock formation and vorticity creation for 3d Euler.
\newblock {\em arXiv:2006.14789v1}, 2020.

\bibitem{Buckmaster2020}
Tristan Buckmaster, Steve Shkoller, and Vlad Vicol.
\newblock Formation of shocks for 2d isentropic compressible Euler.
\newblock  to appear in {\em Communications on Pure and Applied Mathematics}, 2020.

\bibitem{Buckmaster2021}
Tristan Buckmaster, Theodore D. Drivas, Steve Shkoller, and Vlad Vicol.
\newblock Simultaneous development of shocks and cusps for 2D Euler
with azimuthal symmetry from smooth data.
\newblock {\em arXiv:2106.02143v1}, 2021.


\bibitem{Chris2007}
Demetrios Christodoulou.
\newblock {\em The Formation of Shocks in 3-Dimensional Fluids}.
\newblock European Mathematical Society, 2007.

\bibitem{Christodoulou2014}
Demetrios Christodoulou and Shuang Miao.
\newblock {\em Compressible flow and {E}uler's equations}, volume~9 of {\em
  Surveys of Modern Mathematics}.
\newblock International Press, Somerville, MA; Higher Education Press, Beijing,
  2014.


\bibitem{CS96} Constantine G. M. and Savits T. H.\newblock A multivariate Faa di Bruno formula with applications.\newblock {\em Trans. Amer. Math. Soc.} 348 (2), 503–520, 1996.

\bibitem{Deng2003}
Yinbin Deng, Jianlin Xiang, and Tong Yang.
\newblock Blowup phenomena of solutions to {E}uler-{P}oisson equations.
\newblock {\em Journal of Mathematical Analysis and Applications},
  286(1):295--306, 2003.

\bibitem{Eggers2009}
Jens Eggers and Marco~A. Fontelos.
\newblock The role of self-similarity in singularities of partial differential
  equations.
\newblock {\em Nonlinearity}, 22(1):R1--R44, 2009.

\bibitem{Goldreich1980}
P.~Goldreich and S.~V. Weber.
\newblock Homologously collapsing stellar cores.
\newblock {\em The Astrophysical Journal}, 238:991, 1980.

\bibitem{Guo1998}
Yan Guo.
\newblock Smooth irrotational flows in the large to the Euler-Poisson system in
  $\mathbb R^{3+1}$.
\newblock {\em Communications in Mathematical Physics}, 195(2):249--265, 1998.

\bibitem{Guo2021}
Yan Guo, Mahir Had\v{z}i\'{c}, and Juhi Jang.
\newblock Continued gravitational collapse for {N}ewtonian stars.
\newblock {\em Archive for Rational Mechanics and Analysis}, 239(1):431--552, 2021.

\bibitem{Guo2021a}
Yan Guo, Mahir Hadzic, Juhi Jang, and Matthew Schrecker.
\newblock Gravitational collapse for polytropic gaseous stars: Self-similar solutions.
\newblock {\em arXiv:2107.12056 }, 2021.

\bibitem{Guo2016a}
Yan Guo, Lijia Han, and Jingjun Zhang.
\newblock Absence of shocks for one dimensional Euler{\textendash}Poisson
  system.
\newblock {\em Archive for Rational Mechanics and Analysis}, 223(3):1057--1121, 2016.

\bibitem{Guo2014}
Yan Guo, Alexandru~D. Ionescu, and Benoit Pausader.
\newblock Global solutions of certain plasma fluid models in three-dimension.
\newblock {\em Journal of Mathematical Physics}, 55(12):123102, 26, 2014.

\bibitem{Guo2016}
Yan Guo, Alexandru~D. Ionescu, and Benoit Pausader.
\newblock Global solutions of the {E}uler-{M}axwell two-fluid system in 3{D}.
\newblock {\em Annals of Mathematics. Second Series}, 183(2):377--498, 2016.

\bibitem{Guo2011}
Yan Guo and Benoit Pausader.
\newblock Global smooth ion dynamics in the {E}uler-{P}oisson system.
\newblock {\em Communications in Mathematical Physics}, 303(1):89--125, 2011.

\bibitem{Guo1999}
Yan Guo and A.~Shadi Tahvildar-Zadeh.
\newblock Formation of singularities in relativistic fluid dynamics and in
  spherically symmetric plasma dynamics.
\newblock In {\em Nonlinear partial differential equations ({E}vanston, {IL},
  1998)}, volume 238 of {\em Contemp. Math.}, pages 151--161. Amer. Math. Soc.,
  Providence, RI, 1999.

\bibitem{Hadzic2017}
Mahir Had{\v{z}}i{\'{c}} and Juhi Jang.
\newblock Nonlinear stability of expanding star solutions of the radially
  symmetric mass-critical Euler-Poisson system.
\newblock {\em Communications on Pure and Applied Mathematics}, 71(5):827--891, 2017.

\bibitem{Ionescu2013}
Alexandru~D. Ionescu and Benoit Pausader.
\newblock The {E}uler-{P}oisson system in 2{D}: global stability of the
  constant equilibrium solution.
\newblock {\em International Mathematics Research Notices. IMRN}, (4):761--826, 2013.

\bibitem{Jang2008}
Juhi Jang.
\newblock Nonlinear instability in gravitational {E}uler-{P}oisson systems for
  {$\gamma=\frac 65$}.
\newblock {\em Archive for Rational Mechanics and Analysis}, 188(2):265--307, 2008.

\bibitem{Jang2012}
Juhi Jang.
\newblock The two-dimensional {E}uler-{P}oisson system with spherical symmetry.
\newblock {\em Journal of Mathematical Physics}, 53(2):023701, 4, 2012.

\bibitem{Jang2014a}
Juhi Jang.
\newblock Nonlinear instability theory of {L}ane-{E}mden stars.
\newblock {\em Communications on Pure and Applied Mathematics},
  67(9):1418--1465, 2014.

\bibitem{Jang2014}
Juhi Jang, Dong Li, and Xiaoyi Zhang.
\newblock Smooth global solutions for the two-dimensional {E}uler {P}oisson
  system.
\newblock {\em Forum Mathematicum}, 26(3):645--701, 2014.

\bibitem{Li2014}
Dong Li and Yifei Wu.
\newblock The {C}auchy problem for the two dimensional {E}uler-{P}oisson
  system.
\newblock {\em Journal of the European Mathematical Society (JEMS)},
  16(10):2211--2266, 2014.

\bibitem{Li2018}
Hai-Liang Li and Yuexun Wang.
\newblock Formation of singularities of spherically symmetric solutions to the
  3{D} compressible {E}uler equations and {E}uler-{P}oisson equations.
\newblock {\em NoDEA. Nonlinear Differential Equations and Applications},
  25(5):Paper No. 39, 15, 2018.

\bibitem{Liu2019}
Xin Liu.
\newblock On the expanding configurations of viscous radiation gaseous stars:
  the isentropic model.
\newblock {\em Nonlinearity}, 32(8):2975--3011, 2019.

\bibitem{Luk2018}
Jonathan Luk and Jared Speck.
\newblock Shock formation in solutions to the 2{D} compressible {E}uler
  equations in the presence of non-zero vorticity.
\newblock {\em Inventiones Mathematicae}, 214(1):1--169, 2018.

\bibitem{Makino1992}
Tetu Makino.
\newblock Blowing up solutions of the {E}uler-{P}oisson equation for the
  evolution of gaseous stars.
\newblock In {\em Proceedings of the {F}ourth {I}nternational {W}orkshop on
  {M}athematical {A}spects of {F}luid and {P}lasma {D}ynamics ({K}yoto, 1991)},
  volume~21, pages 615--624, 1992.

\bibitem{Makino1990}
Tetu Makino and Beno{\^{\i}}t Perthame.
\newblock Sur les solution {\`{a}} sym{\'{e}}trie sph{\'{e}}rique de l'equation
  d'euler-poisson pour l'evolution d'etoiles gazeuses.
\newblock {\em Japan Journal of Applied Mathematics}, 7(1):165--170, 1990.

\bibitem{Merle2019}
Frank Merle, Pierre Raphael, Igor Rodnianski, and Jeremie Szeftel.
\newblock On smooth self similar solutions to the compressible Euler equations.
\newblock {\em arXiv:1912.10998}, 2019.

\bibitem{Merle2019a}
Frank Merle, Pierre Raphael, Igor Rodnianski, and Jeremie Szeftel.
\newblock On the implosion of a three dimensional compressible fluid.
\newblock {\em arXiv:1912.11009}, 2019.


\bibitem{Rein2003}
Gerhard Rein.
\newblock Non-linear stability of gaseous stars.
\newblock {\em Archive for Rational Mechanics and Analysis}, 168(2):115--130, 2003.

\bibitem{Sideris1985}
Thomas~C. Sideris.
\newblock Formation of singularities in three-dimensional compressible fluids.
\newblock {\em Communications in Mathematical Physics}, 101(4):475--485, 1985.

\bibitem{Wang2014}
Yuexun Wang.
\newblock Formation of singularities to the {E}uler-{P}oisson equations.
\newblock {\em Nonlinear Analysis. Theory, Methods \& Applications. An
  International Multidisciplinary Journal}, 109:136--147, 2014.

\bibitem{Yang2020}
Ruoxuan Yang.
\newblock Shock formation for the Burgers-Hilbert equation.
\newblock {\em arXiv:2006.05568v2}, 2020.

\end{thebibliography}
\end{document}